\documentclass [12pt,a4paper]{article}
\usepackage {amsmath,amsfonts,amssymb,amscd,latexsym,epsfig}
\usepackage {times,mathptm}
\usepackage {theorem}
\usepackage [latin1]{inputenc}
\usepackage [all]{xy}

\input {definitions.tex}

\begin {document}


\thispagestyle {empty}
\markright {INTRODUCTION}

\begin {center}

\vspace* {-14mm}

\textbf {\Large Gromov-Witten invariants of blow-ups}

\vspace {2mm}

{\large Andreas Gathmann}

\vspace {11mm}

\parbox {13cm}{\small
  In the first part of the paper, we give an explicit algorithm to compute
  the (genus zero) \gwinvs of \blowups of an arbitrary convex projective
  variety in some points if one knows the \gwinvs of the original variety. In
  the second part, we specialize to \blowups of $ \PP^r $ and show that many
  invariants of these \blowups can be interpreted as numbers of rational
  curves on $ \PP^r $ having specified global multiplicities or tangent
  directions in the blown-up points. We give various numerical examples,
  including a new easy way to determine the famous multiplicity $ d^{-3} $ for
  $d$-fold coverings of rational curves on the quintic threefold, and, as an
  outlook, two examples of \blowups along subvarieties, whose \gwinvs lead to
  classical multisecant formulas.

  \vspace {7mm}

  \hrule
}

\vspace {3mm}

\end {center}



Over the last few years, \gwinvs of smooth projective varieties have become a
powerful tool in enumerative geometry. Originally applicable only to convex
varieties where the spaces of stable maps have the expected dimension, the
theory is now well-developed for all varieties using virtual fundamental
classes \cite {LT}, \cite {BF}, \cite {B}.

There are at least two motivations to look at \gwinvs of \blowups. Firstly,
a \blowup $ \bX $ of a convex variety $X$ provides an easy example
for a non-convex variety, in the sense that one has reasonably good control
over the stable maps with $ h^1 (C,f^* T_{\bX}) \neq 0 $ since they all must
be such that they intersect the exceptional divisor. Hence this gives a good
class of examples where one can study the effects of virtual fundamental
classes on Gromov-Witten theory. Secondly, curves on the blowup $ \bX $ of a
variety $X$ are closely related to curves on $X$. At least for irreducible
curves not contained in the exceptional divisor, the strict transform of
curves gives a correspondence between curves in $ \bX $ of specified homology
class and curves in $X$ intersecting the blown-up variety with a given
(global) multiplicity. Hence, being able to calculate \gwinvs of blow-ups, one
can hope to solve enumerative problems on $X$ involving multiplicity
conditions at the blown-up variety.

Apart from the last section of this chapter, we will only be concerned with
\blowups of points, since both the calculation and the question of enumerative
significance get very complicated in the case of \blowups of general
subvarieties. Everything will be done over $ \CC $ and for curves of genus
zero.

We first address the question of how one can compute the \gwinvs of \blowups.
For any convex variety $X$, we state and prove an explicit algorithm to
reconstruct all invariants of $ \bX $ from those of $X$ in section \ref
{blowup_algo}. Directly from the algorithm, many of the invariants of $ \bX $
can be seen to vanish or to coincide with others of $X$. This is done in
section \ref {blowup_vanish}. For example, we will show in corollary \ref
{ptexc} that the equality
  \[ \gw [X]{\beta}{\gamma_1 \seq \dots \seq \gamma_n \seq \pt} =
     \gw [\bX]{p^* \beta-E'}{p^* \gamma_1 \seq \dots \seq p^* \gamma_n} \]
holds for $ \beta \in A_1 (X) $ and $ \gamma_i \in A^*(X) $, where $ p: \bX \to
X $ is the \blowup and $ E' $ the class of a line in the exceptional divisor.
As curves in $ \bX $ with homology class $ p^* \beta-E' $ correspond to curves
in $X$ with homology class $ \beta $ intersecting the blown-up point with
multiplicity one, both these invariants are supposed to count curves on $X$
of class $ \beta $ intersecting generic subvarieties representing the $
\gamma_i $ and one additional point in $X$. If the left invariant in fact
counts these curves (which is the case e.g.\ for $ X=\PP^r $ by the
Bertini lemma), then the right invariant also does, and we call this
invariant on $ \bX $ \emph {enumerative} as it has the expected geometric
meaning.

In general, if $ \bX = \bX(s) $ is the \blowup of $X$ at $s$ generic points $
P_1,\dots,P_s $, we will call an invariant on $ \bX $ of the form
  \[ \gw [\bX]{p^* \beta+e_1 E'_1+\dots+e_s E'_s}{
       p^* \gamma_1 \seq \dots \seq p^* \gamma_n} \]
with all $ e_i \le 0 $ enumerative if it counts the number of curves on $X$ of
class $ \beta $ intersecting generic subvarieties representing the $ \gamma_i
$, and in addition passing through each $ P_i $ with global multiplicity $
-e_i $ (see definition \ref {enumdef}). One would then expect these curves to
have $ -e_i $ smooth local branches at every point $ P_i $.

The question whether such a given invariant on $ \bX $ is enumerative or not
is in general very difficult. We will discuss this question in the case
$ \bX = \bP^r(s) $ in sections \ref {blowup_sign} to \ref {blowup_four}.
The results are as follows:
\begin {itemize}
\item If $ s=1 $ then all invariants on $ \bX $ are enumerative. This is shown
  in theorem \ref {enum-1}.
\item If $ r=3 $, $ s \le 4 $, and the invariant contains only point classes
  as incidence conditions, then this invariant is enumerative, except for some
  few cases discussed below. This is shown in theorem \ref {enum-2}.
\item If $ r=3 $ and the invariant contains not only point classes, then it
  is in general not enumerative. This is discussed in section \ref
  {blowup_sign}.
\item If $ r \ge 4 $ and $ s \ge 2 $, then the invariants are ``almost
  never'' enumerative. This is discussed in section \ref {blowup_sign}.
\end {itemize}
In addition, Göttsche and Pandharipande \cite {GP} showed independently that
almost all invariants are enumerative if $ r=2 $. Taking all these results
together, the main point left open is the case $ r=3 $ and $ s \ge 5 $.

In section \ref {blowup_tangency} we show that \gwinvs of \blowups can also be
used to count numbers of curves in $ X=\PP^r $ satisfying certain tangency
conditions: the number of curves in $X$ of class $ \beta $ intersecting
generic representatives of classes $ \gamma_i \in A^*(X) $, and passing in
addition through a given point $ P \in X $ with tangent direction in a given
$k$-codimensional subspace of $ T_{X,P} $ is equal to
\begin {align*}
  \gw [\bX]{p^*\beta-E'}{p^* \gamma_1 \seq \dots \seq p^* \gamma_n \seq
    -(-E)^{k+1}}
    &\qquad \mbox {if $ k<r-1 $}, \\
  \gw [X]{\beta}{\gamma_1 \seq \dots \seq \gamma_n \seq \pt^{\seq 2}}
    - 2 \, \gw [\bX]{p^* \beta-2E'}{p^* \gamma_1 \seq \dots \seq p^* \gamma_n}
    &\qquad \mbox {if $ k=r-1 $},
\end {align*}
see theorem \ref {enum-tang}. Various numerical examples of our results can be
found in section \ref {blowup_examples}. This also includes a very interesting
case of non-enumerative invariants in example \ref {covmult}, namely
  \[ \gw [\bP^3(2)]{d p^* H'-dE'_1-dE'_2}{1} = d^{-3} \]
where $ H' $ is the class of a line in $ \PP^3 $ and the notation $ 1 \in
A^*(X)^{\seq 0} $ means that there are no cohomology classes in the invariant.
This invariant can be shown to coincide with the famous multiplicity with
which multiple coverings get counted in the \gwinvs of the quintic threefold.
Thus our algorithm to compute \gwinvs of \blowups gives a new easy way to
reproduce this result.

We conclude our work with two easy examples of \gwinvs of \blowups of
subvarieties in section \ref {blowup_subv}. In the case of the \blowup
of a space curve $ Y \subset \PP^3 $, we reproduce the well-known (possibly
virtual) number of 3-secants of $Y$ intersecting a fixed line, and the number
of 4-secants of $Y$. In the case of the \blowup of an abelian surface in
$ \PP^4 $, we reproduce the well-known result that the generic abelian surface
in $ \PP^4 $ has 25 6-secants.

This work is part of my PhD thesis written at the University of Hannover.
I would like to thank my advisor Prof.\ K. Hulek for invaluable support and
many helpful discussions. My work has been inspired by my visit of A.
Beauville in Paris, the conference on enumerative geometry in Rome 1997, the
AMS Santa Cruz conference 1995, and in particular by my stay at the
Mittag-Leff\/ler institute last spring during the year on ``Enumerative
geometry and its interactions with theoretical physics''. My work has partly
been financed by the project HCM ERBCHRXCT 940557 (AGE).


\section {Preliminaries} \label {blowup_prelim}

We start by describing the setup and the notation that will be used throughout
the work. For a complex smooth projective variety $X$ of dimension $r$, we
denote by $ \mdf {A_i(X)} $ the algebraic part of $ H_{2i}(X) $ modulo torsion
and by $ \mdf {A^i(X)} $ the algebraic part of $ H^{2i}(X) $ modulo torsion.
These are finitely generated abelian groups. The classes in $ A^i(X) $ will be
said to have \df {codimension} $i$. By abuse of notation, we will often denote
a subvariety of $X$ and its fundamental class in $ A_*(X) $ or $ A^*(X) $ (via
Poincaré duality) by the same symbol if no confusion can result. The
intersection product of two elements $ \gamma $, $ \gamma' $ in $ A^*(X) $ (or
$ A_*(X) $ via Poincaré duality) will be denoted $ \gamma \cdot \gamma' $.
The class of a point will be denoted $ \pt $. If $ X=\PP^r $, the hyperplane
class will be called $ \mdf {H} \in A^1(X) $, and the class of a line will be
called $ \mdf {H'} \in A_1(X) $.

For $ \beta \in A_1 (X) $ an effective homology class and $ n \ge 0 $, we
denote as usual by $ \Mb [n]{X,\beta} $ the moduli spaces of stable maps
of genus zero to $X$ \cite {BM}, and by $ ev_i : \Mb [n]{X,\beta} \to X $ the
evaluation maps. We will sometimes associate to a stable map $ (C,x_1,\dots,
x_n,f) \in \Mb [n]{X,\beta} $ a \df {topology} $ \tau $, by which we mean the
homeomorphism class of the $n$-pointed topological space $ (C,x_1,\dots,x_n) $
together with the data of the homology classes $ f_* [C_i] \in A_1 (X) $ on
each irreducible component $ C_i $ of $C$. This definition can be made much
more precise and formal using the language of graphs \cite {BM}, however then
the notation is likely to get very messy, so we will not make use of it.

These moduli spaces of stable maps possess an expected dimension
  \[ \mdf {\vdim \Mb [n]{X,\beta}} := - K_X \cdot \beta +r+n-3 \]
and a \df {virtual fundamental class} $ \virt {\Mb [n]{X,\beta}} \in A_{\vdim
\Mb [n]{X,\beta}} (\Mb [n]{X,\beta}) $ \cite {LT}, \cite
{BF}, \cite {B}. This class is constructed using the
obstructions $ H^1 (C, f^* T_X) $ for stable maps $ (C,x_1,\dots,x_n,f) \in
\Mb [n]{X,\beta} $. In particular, if these obstructions vanish for all stable
maps in the moduli space, then the virtual fundamental class coincides with the
usual one. There exists a local version of this property too, which follows
immediately from the construction:

\begin {lemma} \label {vfc}
  Let $ (C,x_1,\dots,x_n,f) \in \Mb [n]{X,\beta} $ be a stable map with
  $ h^1 (C,f^* T_X) = 0 $. Then $ (C,x_1,\dots,x_n,f) $ lies in a unique
  irreducible component $Z$ of $ \Mb [n]{X,\beta} $ of dimension $ \vdim
  \Mb [n]{X,\beta} $, and if $R$ denotes the union of all the other
  irreducible components, then
  \begin {gather*}
    \virt {\Mb [n]{X,\beta}} = [Z] + \mbox {some cycle supported on
         $R$}. \tag* {$\Box$}
  \end {gather*}
\end {lemma}

We now come to Gromov-Witten invariants. If $ \gamma_1,\dots,\gamma_n \in
A^*(X) $ are classes on $X$, the associated \gwinv is
  \[ \mdf {\gw [X]{\beta}{\gamma_1 \seq \dots \seq \gamma_n}} :=
       (ev_1^* \gamma_1 \cdot \ldots \cdot ev_n^* \gamma_n) \cdot
       \virt {\Mb [n]{X,\beta}} \in \QQ \]
if $ \sum_{i=1}^n \codim \gamma_i = \vdim \Mb [n]{X,\beta} $, and zero
otherwise.

Concerning the notation, we will often drop the superscript $X$. To shorten
notation, we will often write $ \calT = \gamma_1 \seq \dots \seq \gamma_n $
and call $ \calT \in (A^*(X))^{\seq n} $ a \df {collection of classes}.
Correspondingly, we write $ ev^* \calT $ for $ ev_1^* \gamma_1 \cdot \ldots
\cdot ev_n^* \gamma_n $. If $ X=\PP^r $, the invariant $ \gw {\beta}{\calT} $
is also denoted by $ \gw {d}{\calT} $, where $ \beta = d \, H' $.

We now review briefly the relations among these invariants (see e.g.\ \cite
{FP}), mainly to fix notation for the splitting axiom.

\begin {proposition} \label {axioms}
  Properties of \gwinvs
  \begin {enumerate}
  \item \df {(Mapping to a point)} If $ \beta=0 $, then the invariant is equal
    to the triple intersection product:
      \[ \gw {0}{\gamma_1 \seq \dots \seq \gamma_n} = \begin {cases}
           \gamma_1 \cdot \gamma_2 \cdot \gamma_3 & \mbox {if $n=3$ and
             $ \sum_i \codim \gamma_i = r $,} \\
           0 & \mbox {otherwise.}
         \end {cases} \]
  \item \df {(Fundamental class)} If $ \beta \neq 0 $ and the invariant
    contains the fundamental class of $X$, then the invariant is zero:
      \[ \gw {\beta}{X \seq \calT} = 0 \quad \mbox {for all $ \calT $ and
           all $ \beta \neq 0 $.} \]
  \item \df {(Divisor axiom)} If $ \beta \neq 0 $ and $ \gamma \in A^1(X) $ is
    a divisor, then
      \[ \gw {\beta}{\gamma \seq \calT} = (\gamma \cdot \beta) \,
           \gw {\beta}{\calT} \quad \mbox {for all $ \calT $.} \]
  \item \df {(Splitting axiom)} Choose a homogeneous basis $ \calB = \{ T_0,
    \dots,T_q \} $ of $ A^*(X) $, define $ g = (g_{ij}) $ to be the
    intersection matrix
      \[ \mdf {g_{ij}} = \begin {cases}
           T_i \cdot T_j & \mbox {if $ \codim T_i + \codim T_j = r $,} \\
           0 & \mbox {otherwise,}
         \end {cases} \]
    and let $ g^{-1} = (\mdf {g^{ij}}) $ be the inverse matrix. Choose $ \beta
    \in A_1 (X) $, four classes $ \mu_1,\dots,\mu_4 \in A^* (X) $ and a
    collection $ \calT= \gamma_1 \seq \dots \seq \gamma_n $ of classes such
    that
      \[ \sum_{i=1}^n \codim \gamma_i + \sum_{i=1}^4 \codim \mu_i
           = - K_X \cdot \beta + r + n. \]
    Then we have the equation
    \begin {align*}
      0 = \gw {\beta}{\calT \seq \mu_1 \seq \mu_2 \seq \mu_3 \cdot \mu_4} &
        + \gw {\beta}{\calT \seq \mu_3 \seq \mu_4 \seq \mu_1 \cdot \mu_2} \\
        - \gw {\beta}{\calT \seq \mu_1 \seq \mu_3 \seq \mu_2 \cdot \mu_4} &
        - \gw {\beta}{\calT \seq \mu_2 \seq \mu_4 \seq \mu_1 \cdot \mu_3} \\
        + \sum_{\beta_1,\beta_2 \neq 0} \;
          \sum_{\calT_1,\calT_2} \;
          \sum_{i,j} \;
          g^{ij} \, \Big(
            &\gw {\beta_1}{\calT_1 \seq \mu_1 \seq \mu_2 \seq T_i} \,
            \gw {\beta_2}{\calT_2 \seq \mu_3 \seq \mu_4 \seq T_j} \\
          - &\gw {\beta_1}{\calT_1 \seq \mu_1 \seq \mu_3 \seq T_i} \,
            \gw {\beta_2}{\calT_2 \seq \mu_2 \seq \mu_4 \seq T_j}
          \Big)
    \end {align*}
    where the sum is taken over
    \begin {itemize}
    \item all effective classes $ \beta_1,\beta_2 \in A_1 (X) $ with $ \beta_1
      + \beta_2 = \beta $,
    \item all $ \calT_1 = \gamma_{i_1} \seq \dots \seq \gamma_{i_{n_1}} $ and
      $ \calT_2 = \gamma_{j_1} \seq \dots \seq \gamma_{j_{n_2}} $ such that
      $ i_1 < \dots < i_{n_1} $, $ j_1 < \dots < j_{n_2} $, and $ \{ i_1,\dots,
      i_{n_1} \} \amalg \{ j_1,\dots,j_{n_2} \} = \{ 1,\dots,n \} $ (i.e.\
      ``the classes of $ \calT $ get distributed in all possible ways onto the
      two factors''),
    \item all $ 0 \le i,j \le q $.
    \end {itemize}
    In the sequel we will call this equation $ \mdf {\eqn {\beta}{\calT}{
    \mu_1,\mu_2}{\mu_3,\mu_4}} $.
  \end {enumerate}
\end {proposition}

Now let $ p: \bX=\bX(s) \to X $ be the \blowup of $X$ at $s$ generic points
$ P_1,\dots,P_s \in X $, and let $ E_i $ be the exceptional divisors. Fix a
homogeneous basis $ \mdf {\calB} = \{ T_0,\dots,T_q \} $ of $ A^* (X) $ of
increasing codimension such that $ T_0 = X $ is the fundamental class and $
T_q = \pt $. If we define $ T_{q+1},\dots,T_{\tilde q} $ with $ \tilde q =
q+s(r-1) $ to be the classes
  \[ E_i^k \in A^*(\bX)
       \quad \mbox {where $ 1 \le i \le s, 1 \le k \le r-1 $} \]
(in any order), then
  \[ \mdf {\tilde \calB} =
       \{ p^* T_1,\dots,p^* T_q, T_{q+1},\dots,T_{\tilde q} \} \]
is a homogeneous basis of $ A^*(\bX) $. We call the classes $ p^* T_1,\dots,
p^* T_q $ \df {non-exceptional} and $ T_{q+1},\dots,T_{\tilde q} $ \df
{exceptional}. A collection of classes $ \calT $ will be called non-exceptional
if all its classes are non-exceptional. Since the \gwinvs are multilinear in
the cohomology classes, we will for computational purposes only consider
invariants of the form $ \gw {\beta}{\calT} $ where $ \calT $ is of the form
$ \calT = T_{j_1} \seq \dots \seq T_{j_n}$.

In terms of the basis $ \tilde \calB $, the intersection theory on $ \bX $ is
given by
\begin {align*}
  p^* T_j \cdot p^* T_{j'} &= p^* (T_j \cdot T_{j'}) \\
  p^* T_j \cdot E_i^k &= 0 \\
  E_i^k \cdot E_{i'}^{k'} &= \delta_{i,i'} E_i^{k+k'} \\
  E_i^r &= (-1)^{r-1} \pt
\end {align*}
for $ 1 \le j,j' \le q ;\, 1 \le i,i' \le s ;\, 1 \le k,k' \le r-1 $. If
there is no danger of confusion, we will write the classes $ p^* T_1,\dots,
p^* T_q $ simply as $ T_1,\dots,T_q $.

The homology group $ A_1 (\bX) $ has a canonical decomposition
  \[ A_1(\bX) = A_1(X) \oplus \ZZ \, E'_1 \oplus \dots \oplus \ZZ \, E'_s \]
where $ \mdf{E'_i} $ denotes the class of a line in the exceptional divisor
$ E_i \isom \PP^{r-1} $, such that $ E'_i = -(-E_i)^{r-1} $ via Poincaré
duality. We denote the $ s+1 $ projections onto the summands of the above
decomposition by $ \mdf {d: A_1(\bX) \to A_1(X)} $ and $ \mdf {e_1,\dots,e_s:
A_1(\bX) \to \ZZ} $, and we set $ \mdf {e = e_1 + \dots + e_s} $. If $ X=\PP^r
$, we will identify $ A_1(X) $ with $ \ZZ $ in the obvious way and consider
$d$ as a function $ d: A_1(\bX) \to \ZZ $.

For a homology class $ \beta \in A_1 (\bX) $, we call $ d(\beta) $ the \df
{non-exceptional part} and $ e(\beta) $ the \df {exceptional part}. The class
$ \beta $ is called a \df {non-exceptional class} if $ e_i(\beta) = 0 $ for
all $i$ and a \df {purely exceptional class} if $ d(\beta) = 0 $ and $
e_i(\beta) \neq 0 $ for at least one $i$. For a homology class $ \beta \in A_1
(X) $, we will denote the corresponding non-exceptional class in $ A_1 (\bX) $
also by $ \beta $.

The canonical divisor on $ \bX $ is given by $ K_{\bX} = p^* K_X + (r-1)E $
(see \cite {GH} section 1.4), hence the virtual dimension of the moduli space
$ \Mb [n]{X,\beta} $ is
\begin {align*}
  \vdim \Mb [n]{\bX,\beta} &= -K_{\bX} \cdot \beta +n+r-3 \\
    &= \vdim \Mb [n]{X,d(\beta)} +(r-1) \, e(\beta).
\end {align*}


\section {Calculation of the invariants} \label {blowup_algo}

The aim of this section is to prove the following.

\begin {theorem} \label {reconstruction}
  Let $X$ be a convex variety and $\bX$ the \blowup of $X$ at some points.
  Then there exists an explicit algorithm to compute the \gwinvs of $\bX$
  from those of $X$.
\end {theorem}

The computation is done in three steps. Firstly, we show in lemma \ref
{bXfromX} that all invariants $ \gw [\bX]{\beta}{\calT} $ with $ \beta $ and $
\calT $ non-exceptional are actually equal to the corresponding invariants on
$X$. Secondly, in lemma \ref {purelyexc} we compute the invariants $ \gw
[\bX]{\beta}{\calT} $ with $ \beta $ purely exceptional using a technique
similar to the First Reconstruction Theorem of Kontsevich and Manin. Thirdly,
we state and prove an algorithm that allows one to compute all \gwinvs on $
\bX $ recursively from those obtained in the first two steps.

\begin {lemma} \label {bXfromX}
  Let $ \calT = T_{j_1} \seq \dots \seq T_{j_n} $ be a collection of
  non-exceptional classes and let $ \beta \in A_1 (X) $ be a non-exceptional
  homology class. Then
    \[ \gw [\bX]{\beta}{\calT} = \gw [X]{\beta}{\calT}. \]
  In this case we will say that the invariant $ \gw [\bX]{\beta}{\calT} $ is
  \df {induced by $X$}.
\end {lemma}

\begin {proof}
  Consider the commutative diagram
  \xydiag {
    \Mb [n]{\bX,\beta} \ar[r]^{\phi} \ar[d]_{ev_i}
      & \Mb [n]{X,\beta} \ar[d]_{ev_i} \\
    \bX \ar[r]^p & X
  }
  for $ 1 \le i \le n $. First we show that $ \phi_* \virt {\Mb [n]{\bX,\beta}}
  = \virt {\Mb [n]{X,\beta}} $: since $X$ is convex, $ \Mb [n]{X,\beta} $ is
  a smooth stack of the expected dimension $ d = \vdim \Mb [n]{X,\beta} $. Let
  $ Z_1,\dots,Z_k $ be the connected components of $ \Mb [n]{X,\beta} $, so
  that $ A_d (\Mb [n]{X,\beta}) = \QQ [Z_1] \oplus \dots \oplus \QQ[Z_k] $.
  Since $ \vdim \Mb [n]{\bX,\beta} = d $, we must therefore have
    \[ \phi_* \virt {\Mb [n]{\bX,\beta}} = \alpha_1 [Z_1] + \dots +
                                         \alpha_k [Z_k] \]
  for some $ \alpha_i \in \QQ $.

  To see that all $ \alpha_i = 1 $, pick a stable map $ \cc_i \in Z_i $ whose
  image does not intersect the blown-up points. Then $ \phi^{-1} (\cc_i) $
  consists of exactly one stable map $ \tilde \cc_i $, and the map
  $ \phi: \Mb [n]{\bX,\beta} \to \Mb [n]{X,\beta} $ is a local isomorphism
  around the point $ \tilde \cc_i $. Hence $ \tilde \cc_i $ is a smooth point
  of an irreducible component $ \tilde Z_i $ of $ \Mb [n]{\bX,\beta} $. Denote
  by $ \tilde R_i $ the union of the other irreducible components of $ \Mb
  [n]{\bX,\beta} $. Then, by lemma \ref {vfc},
    \[ \virt {\Mb [n]{\bX,\beta}} = [\tilde Z_i] + \mbox {some cycle supported
         on $ \tilde R_i $}. \]
  Now, since $ \phi: \tilde Z_i \to Z_i $ is a local isomorphism around $
  \tilde \cc_i $, we have $ \phi_* [\tilde Z_i] = [Z_i] $. However, the
  pushforward of a $d$-cycle supported on $ \tilde R_i $ will give no
  contribution to $ \alpha_i $ since $ \cc_i $ and therefore $ Z_i $ is not
  contained in the image of $ \tilde R_i $ under $ \phi $. We conclude that
  all $ \alpha_i = 1 $ and that therefore
  \begin {align*}
    \phi_* \virt {\Mb [n]{\bX,\beta}} &= [Z_1] + \dots + [Z_k] \\
      &= \fund {\Mb [n]{X,\beta}} \\
      &= \virt {\Mb [n]{X,\beta}}.
  \end {align*}
  To complete the proof, note that by the projection formula
  \begin {align*}
    \gw [\bX]{\beta}{\calT}
      &= (\prod_i ev_i^* p^* T_{j_i}) \cdot \virt {\Mb [n]{\bX,\beta}} \\
      &= (\prod_i \phi{}^* ev_i^* T_{j_i}) \cdot \virt {\Mb [n]{\bX,\beta}} \\
      &= (\prod_i ev_i^* T_{j_i}) \cdot \phi_* \virt {\Mb [n]{\bX,\beta}} \\
      &= (\prod_i ev_i^* T_{j_i}) \cdot \virt {\Mb [n]{X,\beta}} \\
      &= \gw [X]{\beta}{\calT}.
  \end {align*}
\end {proof}

\begin {remark} \upshape
  This lemma is actually the only point in the proof of theorem \ref
  {reconstruction} where the convexity of $X$ is needed. Hence, one can
  formulate the theorem also in the following, more general way:

  \begin {sl}
    Let $X$ be a smooth projective variety and $ \bX $ the \blowup of $X$ at
    some points. There exists an explicit algorithm to compute all \gwinvs
    $ \gw [\bX]{\beta}{\calT} $ of $ \bX $ from those where $\beta$ and $
    \calT $ are non-exceptional.
  \end {sl}%

  The proof would be literally the same, just skipping lemma \ref {bXfromX}.
  In fact, it may even be that lemma \ref {bXfromX} also holds for non-convex
  $X$, but I do not know how to prove it in this case.
\end {remark}

\begin {lemma} \label {purelyexc}
  Let $ \calT = T_{j_1} \seq \dots \seq T_{j_n} $ with $ T_{j_i} \in \tilde
  \calB $ be a collection of classes and let $ \beta \in A_1 (\bX) $ be a
  purely exceptional homology class. Then
  \begin {enumerate}
  \item If $ \beta $ is not of the form $ d \cdot E'_i $ for $ d > 0 $ and
    some $ 1 \le i \le s $, then $ \gw [\bX]{\beta}{\calT} = 0 $. Moreover, the
    invariant can only be \nonzero if all classes in $ \calT $ are exceptional
    with support in the exceptional divisor $ E_i $.
  \item $ \gw [\bX]{E'_i}{E_i^{r-1} \seq E_i^{r-1}} = 1 $ for all $ 1 \le i
    \le s $.
  \item All other invariants with purely exceptional homology class can be
    computed recursively.
  \end {enumerate}
\end {lemma}

\begin {proof}
  \begin {enumerate}
  \item This follows easily from the fact that a \gwinv $ \gw [\bX]{\beta}{
    \calT} $ is always zero if there is no stable map in $ \Mb [n]{\bX,\beta}
    $ satisfying the conditions given by $ \calT $.
  \item Note that $ \Mb [2]{\bX,E'_i} \isom \Mb [2]{\PP^{r-1},1} $ and that
    this space is of the expected dimension (which is $ 2r-2 $), hence we do
    not need virtual fundamental classes to compute this invariant. Choose
    two curves $ Y_1,Y_2 \subset X $ intersecting transversally at the blown-up
    point $ P_i $, and let $ \gamma_1,\gamma_2 \in A^{r-1}(X) $ be their
    cohomology classes. Let $ \tilde Y_k $ be the strict transform of $ Y_k $
    for $ k=1,2 $. Then $ \tilde Y_1 $ and $ \tilde Y_2 $ intersect $ E_i $
    transversally at different points, so the invariant
      \[ \gw [\bX]{E'_i}{[\tilde Y_1] \seq [\tilde Y_2]} =
         \gw [\bX]{E'_i}{(\gamma_1 + (-E_i)^{r-1}) \seq (\gamma_2 +
         (-E_i)^{r-1})} \]
    simply counts the number of lines in $ E_i $ through two points in $ E_i
    $, which is 1. Therefore, by the multilinearity of the \gwinvs and by (i)
    we conclude that
    \begin {align*}
      \gw [\bX]{E'_i}{E_i^{r-1} \seq E_i^{r-1}}
        &= \gw [\bX]{E'_i}{(\gamma_1 + (-E_i)^{r-1}) \seq (\gamma_2 +
           (-E_i)^{r-1})} \\
        &= 1.
    \end {align*}
  \item (This is essentially the First Reconstruction Theorem of Kontsevich
    and Manin, see \cite {KM}.) As in (ii) we assume that
    $ \bX = \bP^r(1) $ and that we want to compute the invariant $ \gw
    {d\,E'}{E^{j_1} \seq \dots \seq E^{j_n}} $ for some $d$ and some $ j_i $.
    Consider the equation $ \eqn {d\,E'}{\calT}{E^a,E^b}{E^c,E} $ for some $
    \calT $ consisting of exceptional classes and for some $ 2 \le a \le r-1
    $, $ 2 \le b \le r-1 $, $ 1 \le c \le r-1 $:
    \begin {align*}
      0 &= \gw {d\,E'}{\calT \seq E^a \seq E^b \seq E^c \cdot E} \tag {1} \\
        &+ \gw {d\,E'}{\calT \seq E^c \seq E \seq E^a \cdot E^b} \tag {2} \\
        &- \gw {d\,E'}{\calT \seq E^a \seq E^c \seq E^b \cdot E} \tag {3} \\
        &- \gw {d\,E'}{\calT \seq E^b \seq E \seq E^a \cdot E^c} \tag {4} \\
        &+ \mbox {(terms with homology classes $ d'\,E' $ with $ d' < d $)}.
             \tag {5}
    \end {align*}
    We want to compute the invariants by induction on the degree $d$ and on
    the number of non-divisorial classes in the invariant. Obviously, the
    terms in (5) have lower degree and those in (2) and (4) have same degree
    but a smaller number of non-divisorial classes than (1). The degree of (3)
    is equal to that of (1), and its number of non-divisorial classes is
    not bigger than that of (1). In any case, we can write
    \begin {align*}
      \gw {d\,E'}{\calT \seq E^a \seq E^b \seq E^{c+1}}
        &= \gw {d\,E'}{\calT \seq E^a \seq E^{b+1} \seq E^c} \\
        &+ \mbox {(recursively known terms)}.
    \end {align*}
    Thus if a \gwinv contains at least three non-divisorial classes,
    we can use this equation repeatedly to express $ \gw {d\,E'}{\calT \seq
    E^a \seq E^b \seq E^{c+1}} $ in terms of $ \gw {d\,E'}{\calT \seq
    E^a \seq E^{b+c} \seq E} $ (and recursively known terms), which again has
    fewer non-divisorial classes. This makes the induction work and reduces
    everything to invariants with at most two non-divisorial classes. However,
    since $ \vdim \Mb [n]{\bX,d\,E'} = (r-1)d + r + n - 3 $ and each class
    has codimension at most $r$, it is easy to check that the only such
    invariant is the one calculated in (ii).
  \end {enumerate}
\end {proof}

We now come to the main part of the proof of theorem \ref {reconstruction},
namely the algorithm to compute all invariants on $ \bX $ from those
calculated so far. We will first state the algorithm in such a way that it can
be programmed easily on a computer, and afterwards give the proof that it
really does the job. Many numbers computed using this algorithm can be
found in section \ref {blowup_examples}.

From now on, \gwinvs will always be on $ \bX $ unless otherwise stated, so
we will often write them as $ \gw {\beta}{\calT} $ instead of $ \gw [\bX]{
\beta}{\calT} $.

\begin {algorithm} \label {blowupalgo}
  Suppose one wants to calculate an invariant $ \gw [\bX]{\beta}{\calT} $.
  Assume that the invariant is not induced by $X$ and that $ \beta $ is not
  purely exceptional. We may assume without loss of generality that the sum
  of the codimensions of the non-exceptional classes in $ \calT $ is at least
  $ r+1 $ (hence in particular that there are at least two non-exceptional
  classes) --- otherwise choose a divisor $ \rho \in \calB $ with $ \rho \cdot
  \beta \neq 0 $ (such a $ \rho $ exists because $ \beta $ is not purely
  exceptional) and use $ \calT \seq \rho^{\seq (r+1)} $ instead of $ \calT $,
  which gives essentially the same invariant by the divisor axiom.

  We can further assume without loss of generality that $ \calT $ contains no
  exceptional divisor class and that the classes $ T_{j_1},\dots,T_{j_n} $ in
  $ \calT $ are ordered such that the non-exceptional classes are exactly $
  T_{j_1},\dots,T_{j_m} $, where $ \codim T_{j_1} \ge \dots \ge \codim T_{j_m}
  $. In particular, $ T_{j_1} $ and $ T_{j_2} $ are two non-exceptional
  classes with maximal codimension in $ \calT $.

  We now distinguish the following three cases.
  \begin {itemize}
  \item [(A)] $ n > m $, i.e.\ $ T_{j_n} = E_i^k $ (for some $ 1 \le i \le s $,
    $ 2 \le k \le r-1 $) is an exceptional class. Then use the equation
      \[ \eqn {\beta}{\calT'}{T_{j_1},T_{j_2}}{E_i,E_i^{k-1}} \qquad
           \mbox {where $ \calT' = T_{j_3} \seq \dots \seq T_{j_{n-1}} $}. \]
  \item [(B)] $ n=m $ (i.e.\ there is no exceptional class in $ \calT $),
    $ T_{j_1} = \pt $ and $ \codim T_{j_2} \ge 2 $. Then choose $ \mu, \nu \in
    \calB $ such that $ \codim \mu = 1 $, $ \codim \nu = r-1 $, and $ \mu
    \cdot \nu \neq 0 $. Since the invariant to be computed is not induced by
    $X$, there is an $ i \in \{1,\dots,s\} $ such that $ E_i \cdot \beta \neq
    0 $. Use the equation
      \[ \eqn {\beta}{\calT'}{\mu,\nu}{E_i,T_{j_2}} \qquad
           \mbox {where $ \calT' = T_{j_3} \seq \dots \seq T_{j_n} $}. \]
  \item [(C)] $ n=m $, and it is not true that $ T_{j_1} = \pt $ and $ \codim
    T_{j_2} \ge 2 $. Then again there is an $ i \in \{1,\dots,s\} $ such that
    $ E_i \cdot \beta \neq 0 $. Use the equation
      \[ \eqn {\beta+E'_i}{\calT'}{T_{j_1},T_{j_2}}{E_i,E_i^{r-1}} \qquad
           \mbox {where $ \calT' = T_{j_3} \seq \dots \seq T_{j_n} $}. \]
  \end {itemize}
  Here, ``use equation $ {\calE} $'' means: the \gwinv $ \gw {\beta}{\calT} $
  to be calculated appears in $ {\calE} $ linearly with \nonzero coefficient.
  Solve this equation for $ \gw {\beta}{\calT} $ and compute recursively with
  the same rules all other invariants in this equation that are not already
  known.
\end {algorithm}

\begin {proof} (of theorem \ref {reconstruction})
  Suppose we want to compute an invariant $ \gw {\beta}{\calT} $. If the
  invariant is induced by $X$, it is assumed to be known by lemma \ref
  {bXfromX}. If $ \beta $ is purely exceptional, the invariant is known by
  lemma \ref {purelyexc}. In all other cases, use the algorithm \ref
  {blowupalgo} to compute the invariant recursively. We have to show that the
  equations to be used in fact do contain the desired invariants linearly
  with \nonzero coefficient, and that the recursion stops after a finite
  number of calculations.

  To do this, we will define a partial ordering on pairs $ (\beta,\calT) $
  where $ \beta \in A_1(\bX) $ is an effective homology class and $ \calT $ is
  a collection of cohomology classes. Choose an ordering of the effective
  homology classes in $ A_1 (X) $ such that, for $ \alpha_1,\alpha_2 \neq 0 $
  being two such classes, we have $ \alpha_1 < \alpha_1 + \alpha_2 $ (this is
  possible since the effective classes in $ A_1(X) $ form a semigroup with
  indecomposable zero). For a collection of classes $ \calT = T_{j_1} \seq
  \dots \seq T_{j_n} $, we assume as in the description of the algorithm that
  the classes are ordered such that the non-exceptional classes are exactly
  $ T_{j_1},\dots,T_{j_m} $, where $ \codim T_{j_1} \ge \dots \ge \codim
  T_{j_m} $, and that $ \codim T_{j_1}+\dots+\codim T_{j_m} \ge r+1 $ (by
  possibly adding non-exceptional divisor classes). Then we define
    \[ \mdf {v(\calT)} =
         \min \{ k \;;\; \codim T_{j_1} + \dots + \codim T_{j_k} \ge r+1 \}, \]
  i.e.\ ``the minimal number of non-exceptional classes in $ \calT $ whose
  codimensions sum up to at least $ r+1 $''. With this, we now define the
  partial ordering on pairs $ (\beta,\calT) $ as follows: say that $ (\beta_1,
  \calT_1) < (\beta_2,\calT_2) $ if and only if one of the following holds:
  \begin {itemize}
  \item $ d(\beta_1) < d(\beta_2) $,
  \item $ d(\beta_1) = d(\beta_2) $ and $ v(\calT_1) < v(\calT_2) $,
  \item $ d(\beta_1) = d(\beta_2) $, $ v(\calT_1) = v(\calT_2) $, and $
    e(\beta_1) < e(\beta_2) $.
  \end {itemize}
  Obviously, this defines a partial ordering satisfying the ``descending chain
  condition'', i.e.\ there do not exist infinite chains $ (\beta_1,\calT_1) >
  (\beta_2,\calT_2) > (\beta_3,\calT_3) > \dots $. This means that, to prove
  that the recursion stops after finitely many calculations, it suffices to
  show that the equations in the algorithm compute the desired invariant $ \gw
  {\beta}{\calT} $ entirely in terms of invariants that are either known by
  the lemmas \ref {bXfromX} and \ref {purelyexc} or smaller with respect to
  the above partial ordering. We will do this now for the three cases (A), (B),
  and (C).
  \begin {itemize}
  \item [(A)] The equation reads
    \begin {align*}
      0 &= \gw {\beta}{\calT' \seq T_{j_1} \seq T_{j_2}
             \seq E_i \cdot E_i^{k-1}} \tag {1} \\
        &+ \gw {\beta}{\calT' \seq E_i \seq E_i^{k-1}
             \seq T_{j_1} \cdot T_{j_2}} \tag {2} \\
        &+ \mbox {(no further $ \gw {\beta}{\sth} \, \gw {0}{\sth} $-%
             terms since $ E_i \cdot T_{j_1} = E_i^{k-1} \cdot T_{j_2} = 0 $)}
             \\
        &+ \mbox {(some $ \gw {\beta-d\,E'_i}{\sth} \, \gw {d\,E'_i}{\sth}
             $-terms)} \tag {3} \\
        &+ \mbox {(some $ \gw {\beta_1}{\sth} \, \gw {\beta_2}{\sth}
             $-terms with $ d(\beta_1), d(\beta_2) \neq 0 $)}. \tag {4}
    \end {align*}
    The term (1) is the desired invariant. If the term in (2) is \nonzero, it
    has the same $ d(\beta) $ and smaller $ v(\calT) $, since the two
    non-exceptional classes $ T_{j_1} $, $ T_{j_2} $ of maximal codimensions
    $ \codim T_{j_1} $, $ \codim T_{j_2} $ are replaced by one class of
    codimension $ \codim T_{j_1} + \codim T_{j_2} $. Hence, the term (2) is
    smaller with respect to our partial ordering. The terms in (3) have the
    same $d$, the same or smaller $v$ (note that all non-exceptional classes
    from the original invariant must be in the left invariant $ \gw {\beta-
    d\,E'_i}{\sth} $), and smaller $e$. Finally, the terms in (4) have smaller
    $d$. Hence, all terms in (2), (3) and (4) are smaller with respect to our
    partial ordering.
  \item [(B)] The equation reads
    \begin {align*}
      0 &= \gw {\beta}{\calT' \seq E_i \seq T_{j_2}
             \seq \mu \cdot \nu} \tag {1} \\
        &+ \mbox {(no further $ \gw {\beta}{\sth} \, \gw {0}{\sth} $-%
             terms since $ E_i \cdot T_{j_2} = E_i \cdot \mu = T_{j_2} \cdot
             \nu = 0 $)} \\
        &+ \mbox {(no $ \gw {\beta-d\,E'_i}{\sth} \, \gw {d\,E'_i}{\sth}
             $-terms since $ \gw {d\,E'_i}{\sth} $ would have to contain at
             least} \\
        & \qquad \mbox {one of the non-exceptional classes $ T_{j_2} $,
             $ \mu $, $ \nu $)} \\
        &+ \mbox {(some $ \gw {\beta_1}{\sth} \, \gw {\beta_2}{\sth}
             $-terms with $ d(\beta_1), d(\beta_2) \neq 0 $)}. \tag {2}
    \end {align*}
    Here, obviously, (1) is the desired invariant and the terms in (2) have
    smaller $d$ and are therefore smaller with respect to the partial ordering.
  \item [(C)] The equation reads
    \begin {align*}
      0 &= \gw {\beta+E'_i}{\calT' \seq T_{j_1} \seq T_{j_2}
             \seq \underbrace {E_i \cdot E_i^{r-1}}_{(-1)^{r-1}\pt} }
           \tag {1} \\
        &+ \gw {\beta+E'_i}{\calT' \seq E_i \seq E_i^{r-1}
             \seq T_{j_1} \cdot T_{j_2}} \tag {2} \\
        &+ \mbox {(no further $ \gw {\beta}{\sth} \, \gw {0}{\sth} $-%
             terms)} \\
        &+ \gw {\beta}{\calT' \seq T_{j_1} \seq T_{j_2}
             \seq E_i} \, \underbrace {\gw {E'_i}{E_i \seq E_i^{r-1} \seq
             E_i^{r-1}}}_{=-1} \, (-1)^{r-1} \tag {3} \\
        &+ \mbox {(no further $ \gw {\beta-d\,E'_i}{\sth} \, \gw {d\,E'_i}%
             {\sth} $-terms since there are not enough exceptional} \\
        & \qquad \mbox {classes to put into $ \gw {d\,E'_i}{\sth} $)} \\
        &+ \mbox {(some $ \gw {\beta_1}{\sth} \, \gw {\beta_2}{\sth}
             $-terms with $ d(\beta_1), d(\beta_2) \neq 0 $)}. \tag {4}
    \end {align*}
    Here, (3) is the desired invariant. (4) has smaller $d$, and (2) has the
    same $d$ and smaller $v$, as in case (A)-(2). The term (1) has the same
    $d$, but is not necessarily smaller with respect to the partial ordering.
    We distinguish two cases:
    \begin {enumerate}
    \item If $ \calT' \seq T_{j_1} \seq T_{j_2} $ contains a non-divisorial
      (non-exceptional) class, then the invariant (1) will be computed in the
      next step using rule (B), which expresses it entirely in terms of
      invariants with smaller $d$.
    \item If $ \calT' \seq T_{j_1} \seq T_{j_2} $ contains only divisor
      classes, the invariant (1) will be computed in the next step using (C).
      This time, (2) vanishes (for $ T_{j_1} \cdot T_{j_2} = 0 $ since $
      T_{j_1} = \pt $), (4) has smaller $d$, and (1) will be computed by (B)
      as in (i) in terms of invariants with smaller $d$.
    \end {enumerate}
    Hence, combining (C) with possibly one other application of (B) and/or (C),
    the desired invariant will again be computed in terms of invariants that
    are smaller with respect to the partial ordering.
  \end {itemize}
  This finishes the proof.
\end {proof}

\begin {corollary} \label {proj_rec}
  There exists an explicit algorithm to compute all \gwinvs on $ \bP^r(s) $
  for all $ r \ge 2 $, $ s \ge 1 $.
\end {corollary}

\begin {proof}
  Compute the invariants of $ \PP^r $ using the First Reconstruction Theorem
  \cite {KM}, and then use theorem \ref {reconstruction}.
\end {proof}


\section {A vanishing theorem} \label {blowup_vanish}

We will now prove a vanishing theorem saying that a \gwinv $ \gw {\beta}{\calT
} $ with $ d(\beta) \neq 0 $ and $ e_i(\beta) \ge 0 $ for some $i$ vanishes
under favourable conditions, mainly if $ e_i(\beta) > 0 $ and if there
are ``not too many'' exceptional classes in $ \calT $. The proof of the
proposition is quite involved, but as a reward it is also very sharp in the
sense that numerical calculations on $ \bP^r(1) $ have shown that an invariant
(with non-vanishing $ d(\beta) $ and non-negative $ e(\beta) $) is ``unlikely
to vanish'' if the conditions of the proposition are not satisfied. We will
then apply the proposition to prove corollary \ref {ptexc}, which is a first
hint that \gwinvs on \blowups will lead to enumeratively meaningful numbers.

To state the proposition, we need an auxiliary definition. For $ T \in \tilde
\calB $ and $ 1 \le i \le s $ we define
  \[ \mdf {w_i (T)} = \begin {cases}
       m-1 & \mbox {if $ T=E_i^m $ for some $m$,} \\
       0 & \mbox {otherwise.}
     \end {cases} \]
If $ \calT = T_{j_1} \seq \dots \seq T_{j_n} $ is a collection of classes, we
set $ w_i(\calT) = w_i (T_{j_1}) + \dots + w_i (T_{j_n}) $.

\begin {proposition} \label {vanishing}
  Let $ \beta $ and $ \calT $ be such that for some $ 1 \le i_0 \le s $ the
  following three conditions hold:
  \begin {enumerate}
  \item $ d(\beta) \neq 0 $,
  \item $ w_{i_0} (\calT) > 0 $ or $ e_{i_0} (\beta) > 0 $,
  \item $ w_{i_0} (\calT) < (e_{i_0}(\beta)+1)(r-1) $.
  \end {enumerate}
  Then $ \gw {\beta}{\calT} = 0 $.
\end {proposition}

\begin {proof}
  The proof will be given inductively following the lines of the algorithm
  \ref {blowupalgo}. For invariants induced by $X$ or invariants with
  purely exceptional homology class, the proposition does not say anything, so
  all we have to do is to go through the three equations (A) to (C) and show
  that the statement of the proposition is correct for the invariant to be
  determined if it is correct for all the others.

  For the proof of the proposition, we will refer to the classes $ T_i $ and
  $ T_j $ in the splitting axiom (see proposition \ref {axioms} (iv))
    \[ 0 = \sum \; g^{ij} \, \Big(
         \gw {}{\dots \seq T_i} \, \gw {}{\dots \seq T_j}
         \Big) \]
  as the \df {additional classes} of a certain summand in the equation.

  Assume that we are calculating an invariant $ \gw {\beta}{\calT} $ and that
  a term $ \gw {\beta_1}{\calT_1} \, \gw {\beta_2}{\calT_2} $ occurs in the
  corresponding equation (A), (B), or (C) such that $ (\beta,\calT) $
  satisfies the conditions of the proposition, but neither $ (\beta_1,\calT_1)
  $ nor $ (\beta_2,\calT_2) $ does. We will show that this assumption leads to
  a contradiction.

  We first distinguish the two cases $ w_{i_0} (\calT) > 0 $ and $ e_{i_0}
  (\beta) > 0 $ according to $ (\beta,\calT) $ satisfying (ii).
  \begin {itemize}
  \item \underline {$ w_{i_0} (\calT) > 0 $.} This means that we have an
    exceptional non-divisorial class in the invariant and hence that we are
    in case (A) of the algorithm. Moreover, we can assume that we use case
    (A) of the algorithm with $ i = i_0 $. Since the term in (A)-(2) in the
    proof of theorem \ref {reconstruction} satisfies the conditions of the
    proposition if the desired invariant (A)-(1) does, we only need to
    consider the terms (A)-(3) and (A)-(4).

    From (A)-(1) we know that
      \[ w_i(\calT) = w_i(\calT')+w_i(E_i^k) = w_i(\calT')+k-1, \]
    whereas in all other terms $ \gw {\beta_1}{\calT_1} \, \gw {\beta_2}{
    \calT_2} $ we have
    \begin {gather*}
      w_i (\calT_1) + w_i(\calT_2)
        = w_i(\calT')+w_i(E_i^{k-1}) + \varepsilon (r-2)
        = w_i(\calT')+k-2 + \varepsilon (r-2), \tag {1}
    \end {gather*}
    where $ \varepsilon = 1 $ if the additional classes happen to be
    classes in the exceptional divisor $ E_i $, and $ \varepsilon = 0 $
    otherwise. Combining both equations, we get
      \[ w_i(\calT_1) + w_i(\calT_2) = w_i(\calT) - 1 + \varepsilon (r-2).
           \qquad (*) \]
    Now we again distinguish two cases.
    \begin {itemize}
    \item [(a)] \underline {$ (\beta_1,\calT_1) $ and $ (\beta_2,
      \calT_2) $ satisfy (ii).} If $ (\beta_1,\calT_1) $ does not satisfy (i),
      then $ \beta_1 $ is a purely exceptional class, so all classes in
      $ \calT_1 $ must be exceptional, i.e.\
      \begin {align*}
        w_i(\calT_1) = \vdim \Mb [0]{\bX,\beta_1}
          &= e_i(\beta_1) \, (r-1) + r - 3 \\
          &= (e_i(\beta_1)+1)(r-1) - 2.
      \end {align*}
      So we have the two possibilities
      \begin {align*}
        \mbox {$ (\beta_1,\calT_1) $ does not satisfy (i)} & \;\Rightarrow\;
          w_i(\calT_1) \ge (e_i(\beta_1)+1)(r-1) - 2, \\
        \mbox {$ (\beta_1,\calT_1) $ does not satisfy (iii)} & \;\Rightarrow\;
          w_i(\calT_1) \ge (e_i(\beta_1)+1)(r-1).
      \end {align*}
      The same is true for $ (\beta_2,\calT_2) $. However, since $ \beta $ is
      not purely exceptional, it is not possible that both $ (\beta_1,\calT_1)
      $ and $ (\beta_2,\calT_2) $ do not satisfy (i). We conclude that
      \begin {align*}
        w_i(\calT_1) + w_i(\calT_2)
          &\ge (e_i(\beta_1)+1+e_i(\beta_2)+1)(r-1) - 2 \\
          &= (e_i(\beta)+2)(r-1) - 2 \\
          &> w_i(\calT) + r - 3 \qquad \mbox {since $ (\beta,\calT) $
               satisfies (iii).}
      \end {align*}
      This is a contradiction to (1).
    \item [(b)] \underline {$ (\beta_1,\calT_1) $ does not satisfy (ii),} i.e.\
      $ w_i(\calT_1) = e_i(\beta_1) = 0 $. Since $ w_i(\calT_1) = 0 $, $
      \calT_1 $ does not contain exceptional classes $ E_i^k $ for $ k > 1 $.
      Since $ e_i(\beta_1) = 0 $, $ \calT_1 $ also does not contain $ E_i $
      (otherwise $ \gw {\beta_1}{\calT_1} = 0 $ by the divisor axiom). Hence
      $ \calT_1 $ does not contain $ E_i^k $ for any $k$, and in particular
      we conclude that $ \varepsilon = 0 $ in (1):
      \begin {align*}
        w_i(\calT_2) = w_i(\calT)-1 &< w_i(\calT) \\
          &< (e_i(\beta)+1)(r-1) \\
          &= (e_i(\beta_2)+1)(r-1).
      \end {align*}
      Therefore $ (\beta_2,\calT_2) $ satisfies (iii). It also satisfies (ii),
      since otherwise we would have $ e_i(\beta_1) = e_i(\beta_2) = 0 $ and
      hence get zero by the divisor axiom from the class $ E_i $ in (A). Hence,
      $ (\beta_2,\calT_2) $ cannot satisfy (i), i.e.\ we must be looking at
      the invariants (A)-(3). However, the invariant $ \gw {d'\,E'_i}{\sth} $
      appearing there can never be \nonzero if the additional classes are
      non-exceptional. We reach a contradiction.
    \end {itemize}
  \item \underline {$ e_{i_0}(\beta) > 0 $ and $ w_{i_0} (\calT) = 0 $.} Then
      we can be in any of the cases (A) to (C) of the algorithm. Note that
      $ e_{i_0}(\beta_1) + e_{i_0}(\beta_2) $ is equal to $ e_{i_0}(\beta) $
      or $ e_{i_0}(\beta) + 1 $ (the latter case appearing exactly if we are
      in case (C) and $ i=i_0 $). In any case, it follows that
        \[ e_{i_0}(\beta_1) + e_{i_0}(\beta_2) \ge e_{i_0}(\beta) \ge 1, \]
      hence we can assume without loss of generality that $ e_{i_0}(\beta_1)
      \ge 1 $. In particular, $ (\beta_1,\calT_1) $ satisfies (ii). We are
      going to show that it also satisfies (i) and (iii), which is then a
      contradiction to our assumptions.

      The case that $ (\beta_1,\calT_1) $ does not satisfy (i), i.e.\ that
      $ d(\beta_1) = 0 $, could only occur in (A)-(3) and for $ \beta_1 = d\,
      E'_i $. Since
        \[ 1 \le e_{i_0}(\beta_1) = e_{i_0}(d\,E'_i)
             = d\,\delta_{i,i_0} \]
      we must have $ i=i_0 $. But this means that we have a class $ E_i^k =
      E_{i_0}^k $ in $ \calT $ which is a contradiction to $ w_{i_0}(\calT) =
      0 $. Hence $ (\beta_1,\calT_1) $ must satisfy (i).

      As for (iii), we compute $ w_{i_0}(\calT_1) $. There are no exceptional
      classes $ E_{i_0}^2,\dots,E_{i_0}^{r-1} $ in $ \calT' $ since $ w_{i_0}
      (\calT) = 0 $. Hence the only such classes in $ \calT_1 $ can come from
      \begin {itemize}
      \item the additional classes,
      \item the four special classes used in the equation (A), (B), or (C).
      \end {itemize}
      Both can contribute at most $ r-2 $ to $ w_{i_0}(\calT_1) $, hence
        \[ w_{i_0}(\calT_1) \le 2r-4 < 2(r-1)
             \le (e_{i_0}(\beta_1)+1)(r-1). \]
      Therefore $ (\beta_1,\calT_1) $ also satisfies (iii), arriving at the
      contradiction we were looking for.
  \end {itemize}
\end {proof}

As a corollary we can now prove a relation between the \gwinvs of $ \bX $
that one would expect from geometry. Namely, if we want to express the
condition that curves of homology class $ \beta $ pass through a generic point
in $X$, we expect to be able to do this in two different ways: either we add
the class of a point to $ \calT $, or we blow up the point and count curves
with homology class $ \beta-E' $. The following corollary states that these
two methods will always give the same result, no matter whether the invariants
are actually enumeratively meaningful or not.

\begin {corollary} \label {ptexc}
  Let $ (\beta,\calT) $ be such that, for some $ 1 \le i \le s $, we have
  $ e_i(\beta) = w_i(\calT) = 0 $ and $ d(\beta) \neq 0 $. Then
    \[ \gw {\beta-E'_i}{\calT} = \gw {\beta}{\calT \seq \pt}. \]
\end {corollary}

\begin {proof}
  Consider the equation $ \eqn {\beta}{\calT}{\lambda,\lambda}{E_i,E_i^{r-1}}
  $ for an arbitrary divisor $ \lambda \in \calB $ with $ \lambda \cdot \beta
  \neq 0 $:
  \begin {align*}
    0 &= \gw {\beta}{\calT \seq \lambda \seq \lambda
           \seq E_i \cdot E_i^{r-1}} \tag {1} \\
      &+ \mbox {(no further $ \gw {\beta}{\sth} \, \gw {0}{\sth} $-%
           terms)} \\
      &+ \gw {\beta-E'_i}{\calT \seq \lambda \seq \lambda
           \seq E_i} \, \underbrace {\gw {E'_i}{E_i \seq E_i^{r-1} \seq
           E_i^{r-1}}}_{=-1} \, (-1)^{r-1} \tag {2} \\
      &+ \mbox {(no further $ \gw {\beta-d\,E'_i}{\sth} \, \gw {d\,E'_i}%
           {\sth} $-terms since there are not enough exceptional} \\
      & \qquad \mbox {classes to put into $ \gw {d\,E'_i}{\sth} $)} \\
      &+ \mbox {(some $ \gw {\beta_1}{\sth} \, \gw {\beta_2}{\sth}
           $-terms with $ d(\beta_1), d(\beta_2) \neq 0 $)}. \tag {3}
  \end {align*}
  Using proposition \ref {vanishing}, we will show for any term $ \gw
  {\beta_1}{\calT_1} \, \gw {\beta_2}{\calT_2} $ in (3) that it vanishes.
  Since $ e_i(\beta_1) + e_i(\beta_2) = e_i(\beta) = 0 $, we have without loss
  of generality one of the following cases:
  \begin {itemize}
  \item \underline {$ e_i(\beta_1) = e_i(\beta_2) = 0 $.} Then $ \gw
    {\beta_1}{\calT_1} \, \gw {\beta_2}{\calT_2} = 0 $ by the divisor axiom
    because of the class $ E_i $ in the equation.
  \item \underline {$ e_i(\beta_1) > 0 $.} Then we show that $ (\beta_1,
    \calT_1) $ satisfies conditions (i) to (iii) of the proposition and hence
    vanishes. (i) and (ii) are obvious. As for (iii), the only classes
    contributing to $ w_i(\calT_1) $ can come from
    \begin {itemize}
    \item the additional classes,
    \item the special class $ E_i^{r-1} $ used in the equation.
    \end {itemize}
    Both can contribute at most $ r-2 $ to $ w_i(\calT_1) $, hence
    \[ w_i(\calT_1) \le 2r-4 < 2(r-1) \le (e_i(\beta_1)+1)(r-1). \]
    Therefore $ (\beta_1,\calT_1) $ also satisfies (iii).
  \end {itemize}
  Now that we know that all terms in (3) vanish, the above equation becomes
    \[ \gw {\beta}{\calT \seq \lambda \seq \lambda \seq E_i \cdot E_i^{r-1}}
       = \gw {\beta-E'_i}{\calT \seq \lambda \seq \lambda
           \seq E_i} \, (-1)^{r-1}. \]
  Since $ E_i \cdot E_i^{r-1} = (-1)^{r-1} \pt $ and $ E_i \cdot (\beta-E'_i)
  =1 $, the corollary follows.
\end {proof}


\section {Enumerative significance --- general remarks} \label {blowup_sign}

After having computed all \gwinvs on \blowups of projective space (see
corollary \ref {proj_rec}), we now come to the question of enumerative
significance of the invariants. For most of the time, we will be concerned with
invariants $ \gw [\bX]{\beta}{\calT} $ where $ \calT $ is non-exceptional,
leading to numbers of curves on $X$ intersecting the blown-up points with
prescribed multiplicities. Only in section \ref {blowup_tangency} we will
consider some invariants containing exceptional classes in $ \calT $, leading
to numbers of curves on $X$ with certain tangency conditions.

For the rest of the chapter, we will only work with $ \bX=\bP^r(s) $. We start
by giving a precise definition of an enumeratively significant invariant.

\begin {definition} \label {enumdef}
  Let $ \beta \in A_1 (\bX) $ a homology class with $ d(\beta) \neq 0 $ and
  $ e_i(\beta) \le 0 $, and let $ \calT = \gamma_1 \seq \dots \seq \gamma_n
  $ be a collection of non-exceptional effective classes $ \gamma_i \in A^{\ge
  1} (X) $ such that $ \sum_i \codim \gamma_i = \vdim {\Mb [n]{\bX,\beta}} $.

  Then we call the \gwinv $ \gw [\bX]{\beta}{\calT} $ \df {enumerative} if,
  for generic subschemes $ V_i \subset \bX $ with $ [V_i] = \gamma_i $, it is
  equal to the number of irreducible stable maps $ (C,x_1,\dots,x_n,f) $ with
  $f$ being generically injective, $ f_* [C] = \beta $, and $ f(x_i) \in V_i $
  for all $i$ (where each such stable map is counted with multiplicity one).
\end {definition}

Note that irreducible stable maps $ (C,x_1,\dots,x_n,f) $ on $ \bX $ of
homology class $ \beta $ with $f$ generically injective correspond bijectively
to irreducible curves in $ \bX $ of homology class $ \beta $, and hence via
strict transform to irreducible curves in $X$ of homology class $ d(\beta) $
intersecting the blown-up points $ P_i $ with global multiplicities $ -e_i
(\beta) $. Hence it is clear that we can also give the following
interpretation of enumerative invariants:

\begin {lemma}
  If $ \gw {\beta}{\calT} $ is enumerative, then for generic subschemes $ V_i
  \subset \bX $ with $ [V_i] = \gamma_i $, it is equal to the number of
  irreducible rational curves $ C \subset X $ of homology class $ d(\beta) $
  intersecting all $ V_i $, and in addition passing through each $ P_i $
  with global multiplicity $ -e_i(\beta) $. Every such curve is counted with
  multiplicity $ \sharp (C \cap V_1) \cdot \ldots \cdot \sharp (C \cap V_n) $.
\end {lemma}

In general, one would then expect these curves to have $ -e_i $ smooth local
branches at every point $ P_i $.

We will now give an overview of the results about enumerative significance
of \gwinvs on $ \bP^r(s) $. Assume that $ d(\beta) \neq 0 $, $ e_i(\beta) \le
0 $, and that $ \calT $ is a collection of non-exceptional effective classes.

\begin {enumerate}
  \item If $ s=1 $ then $ \gw {\beta}{\calT} $ is enumerative. This will
    be shown in theorem \ref {enum-1}.
  \item If $ r=2 $ then $ \gw {\beta}{\calT} $ is enumerative if $ e_i(\beta)
    \in \{-1,-2\} $ for some $i$ or $ \calT $ contains at least one point
    class. This has been proven by L. Göttsche and R. Pandharipande in \cite
    {GP}.
  \item If $ r=3 $, $ s \le 4 $, and $ \calT $ contains only point classes,
    then $ \gw {\beta}{\calT} $ is enumerative if and only if $ \beta $ is not
    equal to $ d\,H'-d\,E'_i-d\,E'_j $ for some $ d \ge 2 $ and $ i \neq j $
    with $ 1 \le i,j \le s $. We will prove this in theorem \ref {enum-2}.
  \item If $ r=3 $ and $ \calT $ contains not only point classes, then $ \gw
    {\beta}{\calT} $ is in general not enumerative.
  \item If $ r \ge 4 $ and $ s \ge 2 $ then $ \gw {\beta}{\calT} $ is
    ``almost never'' enumerative.
\end {enumerate}

We start our study of enumerative significance by showing the origin of
potential problems with enumerative significance, thereby giving
counterexamples to enumerative significance in the cases (iv) and (v) above.

The most obvious problem is that a stable map $ (C,x_1,\dots,x_n,f) $ may be
reducible, with some of the components mapped to the exceptional divisor.
The part of the moduli space corresponding to such curves will in general have
too big dimension. For example, consider the case $ \bX=\bP^3(1) $, $ \beta =
4H' $. Stable maps in $ \Mm [0]{\bX,\beta} $ will not intersect the
exceptional divisor at all, hence $ \Mm [0]{\bX,\beta} $ has the expected
dimension. However, consider reducible curves $ C = C_1 \cup C_2 $ where $f$
is of homology class $ 4H'-3E' $ on $ C_1 $ and of homology class $ 3E' $ on
$ C_2 $. These can be depicted as follows:
\begin {center} \unitlength 1cm \begin {picture}(4,4)
  \put (0,0){\epsfig {file=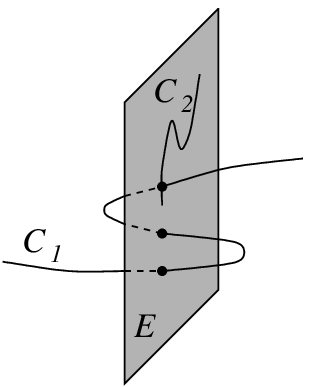,width=3cm}}
\end {picture} \end {center}
The space of such curves $ C_1 $ is (at least) of dimension $ \vdim \Mb [0]{
\bX,4H'-3E'} = 4 \cdot 4 - 3 \cdot 2 = 10 $, the space of curves $ C_2 $ of
homology class $ 3E' $ through a given point (namely one of the points of
intersection of $ C_1 $ with $E$) is of dimension $ 3 \cdot 3 - 1 - 1 = 7 $
(note that $ E \isom \PP^2 $). Hence the part of the moduli space $ \Mb [0]{
\bX,\beta} $ corresponding to those curves has dimension (at least) 17, but we
have $ \vdim \Mb [0]{\bX,\beta} = 4 \cdot 4 = 16 $. Note that this is in
agreement with the fact that these curves certainly cannot be deformed into
smooth quartics not intersecting the exceptional divisor, hence they are not
contained in the closure of $ \Mm [0]{\bX,\beta} $ in $ \Mb [0]{\bX,\beta} $.

However, this will cause no problems when computing \gwinvs, since, intuitively
speaking, the curve $ C_2 $ cannot satisfy any incidence conditions with
generic non-exceptional varieties. So if we try to impose $ \vdim \Mb
[0]{\bX,\beta} = 16 $ non-exceptional conditions on these curves, we will get
zero, since the curve $ C_1 $ can satisfy at most 10 of the conditions and
$ C_2 $ can satisfy none at all. For a mathematically more precise statement
of this fact, see proposition \ref {dim-im} (i) which is the important step in
the proof of enumerative significance in the case of only one \blowup.

When we consider more than one \blowup, things get more complicated, since
then for example multiple coverings of the lines joining the blown-up points
will cause problems. As an example, consider $ \bX = \bP^r(2) $, $ \beta =
(d+q)\,H'-q\,E'_1-q\,E'_2 $ for some $ r \ge 2 $, $ d \ge 1 $, $ q \ge 2 $, and
look at reducible stable maps as above with $ C_1 $ of homology class $ d\,H'
$ and $ C_2 $ of homology class $ q\,H'-q\,E'_1-q\,E'_2 $, being a $q$-fold
covering of the strict transform of the line between $ P_1 $ and $ P_2 $:
\begin {center} \unitlength 1mm \begin {picture}(48,24)
  \put (0,0){\epsfig {file=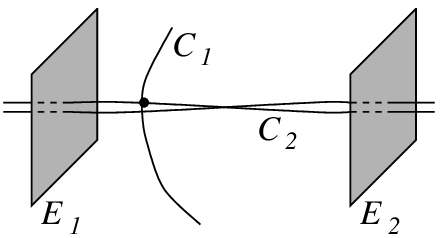,width=48mm}}
\end {picture} \end {center}
We have just learned that $ C_2 $ for itself will make no problems, since no
generic (non-divisorial) non-exceptional incidence conditions can be satisfied
on this component. However, it may well happen that the dimension of the
moduli space of curves $ C_1 $ meeting the line through $ P_1 $ and $ P_2 $
(i.e.\ $ \vdim \Mb [0]{\bX,d\,H'}-(r-2) $) is \emph {bigger} than that of both
components together:
\begin {align*}
  \vdim \Mb [0]{\bX,d\,H'}-(r-2) &= (r+1)d+r-3-(r-2), \\
  \vdim \Mb [0]{\bX,\beta} &= (r+1)d+(1-q)(r-3), \\
  \Rightarrow \vdim \Mb [0]{\bX,d\,H'} - (r-2) - \vdim \Mb [0]{\bX,\beta}
    &= \underline {(q-1)(r-3)-1}.
\end {align*}
If this last number is non-negative, we will obviously get non-wanted
contributions to our \gwinvs from these reducible curves, since all $ \vdim
\Mb [0]{X,\beta} $ conditions that we impose on the curve can be satisfied on
$ C_1 $. This will always happen if $ r \ge 4 $, showing that in this case
there is no chance of getting enumerative invariants. The reader who wants to
convince himself of this fact numerically can find some obviously
non-enumerative invariants of this kind in example \ref {p4-1}. For $ r=3 $,
we will see that multiple coverings of lines joining blown-up points only make
problems if they form the only component of an irreducible curve, see
theorem \ref {enum-2} and example \ref {p3-2}. In fact, in the case where
$ \beta = d\,H'-d\,E'_1-d\,E'_2 $, such that we ``count'' $d$-fold coverings of
lines, we get other important invariants, see example \ref {covmult}.

Since the case of $ \bP^4(s) $ for $ s \ge 2 $ will not lead to enumerative
invariants and the case of $ \bP^2(s) $ has been studied almost exhaustively
in \cite {GP}, it only remains to look at \blowups of $ \PP^3 $. We will look
at the case $ \bX = \bP^3(4) $ in detail in section \ref {blowup_four} (which
then includes, of course, also the cases $ \bX = \bP^3(s) $ with $ s<4 $).
Here, in analogy to the situation discussed above, one gets problems with too
big dimensions for reducible curves as above, where $ C_2 $ is now a curve
contained in a plane spanned by three of the blown-up points. These problems
arise in particular because in this case it is no longer true that $ C_2 $ can
satisfy no incidence conditions. To be more precise, $ C_2 $ can satisfy
incidence conditions with generic curves, but \emph {not} with generic points
in $ \bP^3(4) $. This is the reason why we have to make the assumption
that all cohomology classes in the invariant are point classes (see
theorem \ref {enum-2}). If we do not assume this, we can again easily get
non-enumerative invariants, e.g.\ $ \gw [\bP^3(4)]{4H'-2E'_1-2E'_2-2E'_3}{
(H^2)^{\seq 4}} = -1 $, to mention the easiest one.

In the remainder of this section, we will prove some statements about
irreducible curves in \blowups that will be needed for both cases $ \bP^r(1) $
and $ \bP^3(4) $. We start by computing $ h^1 (\PP^1,f^* T_{\bX}) $ in the
next two lemmas.

\begin {lemma} \label {cohbX}
  Let $ p: \bX \to X $ be the \blowup of a smooth variety at some points
  $ P_1,\dots,P_s $ and let $ E = E_1 \cup \dots \cup E_s $ be the exceptional
  divisor. Let $C$ be a smooth curve and $ f: C \to \bX $ a map such that
  $ f(C) \not\subset E $. Then there is an injective morphism of sheaves on
  $ \bX $
    \[ f^* p^* T_X (-f^* E) \to f^* T_{\bX} \]
  which is an isomorphism away from $ f^{-1} (E) $.
\end {lemma}

\begin {proof}
  Since $ E = \{ P_1,\dots,P_s \} \times_X \bX $, we have $ i^* \Omega_{
  \bX/X} = \Omega_{E/\{ P_1,\dots,P_s \}} = \Omega_E $ where $ i: E \to \bX $
  is the inclusion. As $ \Omega_{\bX/X} $ has support on $E$, this can be
  rewritten as $ i_* \Omega_E = \Omega_{\bX/X} $. Hence, there is an exact
  sequence of sheaves on $\bX$
    \[ 0 \to p^* \Omega_X \to \Omega_{\bX} \to i_* \Omega_E \to 0. \]
  Dualizing, we get
    \[ 0 \to T_{\bX} \to p^* T_X \to \EExt^1 (i_* \Omega_E, \OO_X) \to 0. \]
  By duality (see \cite {HR} theorem III 6.7), we have
    \[ \EExt^1 (i_* \Omega_E, \OO_X) = i_* \EExt^1 (\Omega_E,N_{E/\bX})
         = i_* T_E (-1) \]
  where $ \OO(-1) := \OO_{E_1}(-1) \seq \dots \seq \OO_{E_s}(-1) $.
  Therefore we get a morphism $ p^* T_X \to i_* T_E(-1) $ which we can
  restrict to $E$ to get a morphism $ p^* T_X|_E \to i_* T_E(-1) $ fitting
  into a commutative diagram
  \[ \begin {CD}
    0 @>>> p^* T_X (-E) @>>> p^* T_X @>>>  p^* T_X|_E  @>>> 0  \\
    @.           @.             @|            @VVV          @. \\
    0 @>>>   T_{\bX}    @>>> p^* T_X @>>> i_* T_E (-1) @>>> 0.
  \end {CD} \]
  From this we can deduce the existence of an injective map $ p^* T_X(-E) \to
  T_{\bX} $ which is clearly an isomorphism away from $E$. Applying the functor $ f^*
  $ we get the desired morphism $ f^* p^* T_X (-f^* E) \to f^* T_{\bX} $. Since
  the image of $f$ is not contained in $E$, this morphism is also injective
  and an isomorphism away from $ f^{-1}(E) $.
\end {proof}

\begin {lemma} \label {coh-irr}
  Let $ C = \PP^1 $, $ \bX = \PP^r(s) $, $ f: C \to \tilde X $
  a morphism, $ \beta = f_*[C] \in A_1(\bX) $, and $ \varepsilon \in \{0,1\} $.
  \begin {enumerate}
  \item If $ f(C) \not\subset E $ or $f$ is a constant map then $ h^1
    (C,f^* T_{\bX}(-\varepsilon)) = 0 $ whenever $ d(\beta)+e(\beta) \ge 0 $.
    (Here, $ f^* T_{\bX}(-\varepsilon) $ is to be interpreted as $ f^* T_{\bX}
    \seq \OO_C (-\varepsilon) $.) In particular, this always holds for $
    s=1 $ (since then $ d(\beta)+e(\beta) = \deg f^* (H-E) $ and $ f^*(H-E) $
    is an effective divisor on $C$).
  \item If $ f(C) \subset E $ and the map $ f: C \to E \cong \PP^{r-1} $ has
    degree $ e>0 $ then
      \[ h^1 (C,f^* T_{\tilde X}(-\varepsilon)) = e + \varepsilon - 1. \]
  \end {enumerate}
\end {lemma}

\begin {proof}
  \begin {enumerate}
  \item If $f$ is a constant map then the assertion is trivial, so assume that
    $ f(C) \not\subset E $ and set $ d=\deg f^* H $, $ e=- \deg f^* E $. By
    lemma \ref {cohbX} we have an exact sequence
      \[ 0 \to f^* p^* T_X (e) \to f^* T_{\bX} \to Q \to 0 \]
    with some sheaf $Q$ on $C$ with zero-dimensional support. Hence to prove
    the lemma it suffices to show that $ h^1 (C,f^* p^* T_X (e-\varepsilon)) =
    0 $. But this follows from the Euler sequence on $ \PP^r $ pulled back to
    $C$ and twisted by $ \OO_C (e-\varepsilon) $:
          \[ 0 \to \OO_C (e-\varepsilon)
           \to (r+1) \, \OO_C (d+e-\varepsilon)
           \to f^* p^* T_X (e-\varepsilon)
           \to 0 \]
    since $ d+e-\epsilon \ge -1 $ by assumption.
  \item We consider the normal sequence
      \[ 0 \to T_E \to i^* T_{\bX} \to N_{E/\bX} \to 0. \]
    As $ N_{E/\bX} = \OO_E (-1) $, pulling back to $C$ and twisting by
    $ \OO_C (-\varepsilon) $ yields
    \begin {gather*}
      0 \to f^* T_E (-\varepsilon)
        \to f^* T_{\bX} (-\varepsilon)
        \to \OO_C (-e-\varepsilon)
        \to 0. \tag {1}
    \end {gather*}
    In complete analogy to (i), it follows by the Euler sequence of $ E \isom
    \PP^{r-1} $
      \[ 0 \to \OO_C (-\varepsilon)
           \to r \, \OO_C (e-\varepsilon)
           \to f^* T_E (-\varepsilon)
           \to 0 \]
    that $ h^1 (C, f^* T_E (-\varepsilon)) = 0 $. Hence we deduce from (1)
    that
      \[ h^1 (f^* T_{\bX} (-\varepsilon)) = h^1 (C,\OO_C (-e-\varepsilon))
           = e + \varepsilon - 1. \]
  \end {enumerate}
\end {proof}

We now come to the Bertini lemma \ref {bertini} which is our main tool to
prove the transversality of the intersection products in the \gwinvs.

\begin {lemma} \label {bertini-0}
  Let $M$ be a scheme of finite type and $ f: M \to \PP^r $ a morphism. Then,
  for a generic hyperplane $ H \subset \PP^r $, we have:
  \begin {enumerate}
    \item $ f^{-1}(H) $ is (empty or) of pure codimension 1 in $M$.
    \item If $M$ is smooth then the divisor $ f^{-1}(H) $ is a smooth
      subscheme of $M$ counted with multiplicity one.
  \end {enumerate}
\end {lemma}

\begin {proof}
  See e.g.\ \cite {J} corollary 6.11.
\end {proof}

\begin {lemma} \label {bertini-1}
  Let $M$ be a scheme of finite type, $X$ a smooth, connected, projective
  scheme, and $ f: M \to X $ a morphism. Let $L$ be a base point free linear
  system on $X$. Then, for generic $ D \in L $, we have:
  \begin {enumerate}
    \item $ f^{-1}(D) $ is (empty or) purely 1--codimensional.
    \item If $M$ is smooth then the divisor $ f^{-1}(D) $ is a smooth
      subscheme of $M$ counted with multiplicity one.
  \end {enumerate}
\end {lemma}

\begin {proof}
  The base point free linear system $L$ on $X$ gives rise to a morphism
  $ s: X \to \PP^m $ where $ m = \dim L $. Composing with $f$ yields a
  morphism $ M \to \PP^m $, and the divisors $ D \in L $ correspond to the
  inverse images under $s$ of the hyperplanes in $ \PP^m $. Hence, the
  statement follows from lemma \ref {bertini-0}, applied to the map $ M \to
  \PP^m $.
\end {proof}

\begin {lemma} \label {bertini}
  Let $M$ be a Deligne-Mumford stack of finite type, $X$ a smooth, connected,
  projective scheme and $ f_i : M \to X $ morphisms for $ i = 1,\dots,n $. Let
  $ \gamma_i \in A^{c_i}(X) $ be cycles of codimensions $ c_i \ge 1 $ on $X$
  that can be written as intersection products of divisors on $X$
    \[ \gamma_i = [D'_{i,1}] \cdot \dots \cdot [D'_{i,c_i}]
         \qquad \mbox {($ i = 1,\dots,n $)} \]
  such that the complete linear systems $ | D'_{i,j} | $ are base point free
  (this always applies, for example, for effective cycles in the case $
  X=\PP^r $). Let $ c = c_1 + \cdots + c_n $. Then, for generic $ D_{i,j}
  \in | D'_{i,j} | $, we have:
  \begin {enumerate}
    \item $ V_i := D_{i,1} \cap \cdots \cap D_{i,c_i} $ is smooth of
      pure codimension $ c_i $ in $X$, and the intersection is transverse.
      In particular, $ [V_i] = \gamma_i $.
    \item $ V := f_1^{-1}(V_1) \cap \cdots \cap f_n^{-1}(V_n) $ is of pure
      codimension $c$ in $M$. In particular, if $ \dim M < c $ then $ V =
      \emptyset $.
    \item If $ \dim M = c $ and $M$ contains a dense, open, smooth substack
      $U$ such that each geometric point of $U$ has no non-trivial
      automorphisms then $V$ consists of exactly $ (f_1^*\gamma_1 \cdot \ldots
      \cdot f_n^*\gamma_n)[X] $ points of $M$ which lie in $U$ and are counted
      with multiplicity one.
  \end {enumerate}
\end {lemma}

\begin {proof}
  \begin {enumerate}
  \item follows immediately by recursive application of lemma \ref {bertini-0}
    to the scheme $X$.
  \item If $M$ is a scheme, then the statement follows by recursive
    application of lemma \ref {bertini-1}. If $M$ is a Deligne-Mumford stack,
    then it has an étale cover $ S \to M $ by a scheme $S$, so (ii) holds for
    the composed maps $ S \to M \to X $. But since the map $ S \to M $ is
    étale, the statement is also true for the maps $ M \to X $.
  \item A Deligne-Mumford stack $U$ whose generic geometric point has no
    non-trivial automorphisms always has a dense open substack $U'$ which is a
    scheme (see e.g.\ \cite {V}. To be more precise, $U$ is a functor and
    hence an algebraic space (\cite {DM} ex.\ 4.9), but an algebraic space
    always contains a dense open subset $U'$ which is a scheme (\cite {Kn}
    p.\ 25)). Since $U'$ is dense in $M$ and therefore $ M \backslash U' $ has
    smaller dimension, applying (ii) to the restrictions $ f_i |_{M \backslash
    U'} : M \backslash U' \to X $ gives that $V$ is contained in the smooth
    scheme $U'$, hence it suffices to consider the restrictions $ f_i |_{U'} :
    U' \to X $. But since $U'$ is a smooth scheme, we can apply lemma \ref
    {bertini-1} (ii) recursively and get the desired result.
  \end {enumerate}
\end {proof}

As we needed for lemma \ref {bertini} (iii) that the generic element of $M$
has no non-trivial automorphisms, we now give a criterion under which
circumstances this is satisfied for our moduli spaces of stable maps.

\begin {lemma} \label {noauto}
  Let $ \bX=\PP^r(s) $ and $ \beta \in A_1 (\bX) $ with $ d(\beta) > 0 $ and $
  d(\beta)+e(\beta) \ge 0 $. Assume that $ \beta $ is not of the form $ d\,H'
  -d\,E'_i $ for $ 1 \le i \le s $ and $ d \ge 2 $. Then, if $ \Mm [n]{\bX,
  \beta} $ is not empty, it is a smooth stack of the expected dimension, and
  if $ \cc=(C,x_1,\dots,x_n,f) $ is a generic element of $ \Mm [n]{\bX,\beta}
  $ then $\cc$ has no automorphisms and $f$ is generically injective.
\end {lemma}

\begin {proof}
  Set $ d=d(\beta) $ and $ e=e(\beta) $. We can assume that $ e \le 0 $ since
  otherwise $ \Mm [n]{\bX,\beta} $ is empty.

  It follows from lemma \ref {coh-irr} (i) that $ \Mm [n]{\bX,\beta} $ is a
  smooth stack of the expected dimension. Note that an irreducible stable map
  can only have automorphisms if it is a multiple covering map onto its image.
  Therefore it suffices if we compute, for all $ N \ge 2 $, the dimension of
  the subspace $ Z_N \subset \Mm [n]{\bX,\beta} $ consisting of $N$-fold
  coverings and show that it is smaller than the dimension of $ \Mm [n]{\bX,
  \beta} $.

  So assume that $ N \ge 2 $ and that $ Z_N \neq \emptyset $, so that $
  \beta = N \beta' $ for some $ \beta' \in A_1 (\bX) $. We set $ d' = d(\beta'
  ) $ and $ e' = e(\beta') $. Since $ d'+e' \ge 0 $, we can apply lemma \ref
  {coh-irr} (i) to see that the space of stable maps of homology class $
  \beta' $ is of the expected dimension $ (r+1) d' + (r-1) e' + r+n-3 $. The
  dimension of $ Z_N $ is exactly bigger by $ 2N-2 $ because of the moduli of
  the covering. Hence we have
  \begin {align*}
    \dim Z_N &= (r+1) d' + (r-1) e' + r+n-3 + 2N-2 \\
      &= (r+1)d + (r-1)e + r+n-3 + ((r+1)d' + (r-1)e') (1-N) + 2N-2 \\
      &= \dim \Mm [n]{\bX,\beta} + ((r+1)d' + (r-1)e' - 2) (1-N).
  \end {align*}
  Therefore, to prove the lemma, it suffices to show that $ (r+1)d' + (r-1)e'
  > 2 $. We distinguish two cases:
  \begin {itemize}
  \item If $ e'=0 $, then
       \[ (r+1)d' + (r-1)e' = (r+1)d' \ge (2+1) \cdot 1 = 3 > 2. \]
  \item If $ e' \le -1 $, then
      \[ (r+1)d' + (r-1)e' = (r+1)(d'+e') - 2e' \ge -2e' \ge 2, \]
    but if we had equality, this would mean $ d'+e'=0 $ and $ e'=-1 $, hence
    $ \beta' = H'-E'_i $ for some $i$ and therefore $ \beta = N\,H'-N\,E'_i $,
    which is the case we excluded in the lemma.
  \end {itemize}
  This finishes the proof.
\end {proof}


\section {Enumerative significance --- the case $ \bP^{\symbol {114}}(1) $}
  \label {blowup_one}

In this section we will prove that all invariants $ \gw {\beta}{\calT} $
on $ \bX = \bP^r(1) $ are enumerative. We start with the computation of $ h^1
(C,f^* T_{\bX}) $ for arbitrary stable maps. To state the result, we need the
following definition: for any prestable map $ (C,x_1,\dots,x_n,f) $ to $ \bX $
we define $ \eta (C,f) $ to be ``the sum of the exceptional degrees of all
irreducible components of $C$ which are mapped into $E$'', i.e.\
  \[ \mdf {\eta (C,f)} := \sum_{C'} \;\, \{ \mbox {
       $ e \; | \; C' $ is an irreducible component of $C$
       such that $ f_* [C'] = e \, E' $
     }\}. \]
Obviously, $ \eta (C,f) $ only depends on the topology $ \tau $ of the
prestable map in the sense of section \ref {blowup_prelim}, so we will write
$ \eta (\tau)=\eta(C,f) $.

\begin {lemma} \label {coh-red}
  Let $C$ be a prestable curve, $ \bX = \bP^r(1) $, and $ f: C \to \bX $ a
  morphism. Then $ h^1 (C,f^* T_{\bX}) \le \eta (C,f) $, with strict
  inequality holding if $ \eta (C,f) > 0 $.
\end {lemma}

\begin {proof}
  The proof is by induction on the number of irreducible components of $C$. If
  $C$ itself is irreducible, the statement follows immediately from lemma
  \ref {coh-irr} for $ \varepsilon=0 $.

  Now let $C$ be reducible, so assume $ C = C_0 \cup C' $ where $ C' \cong
  \PP^1 $, $ C_0 \cap C' = \{Q\} $, and where $ C_0 $ is a prestable curve
  for which the induction hypothesis holds. If $ \eta (C,f) > 0 $, we can
  arrange this such that $ \eta (C_0,f_0) > 0 $.

  Consider the exact sequences
  \begin {gather*}
    0 \to f^* T_{\tilde X}
      \to f_0^* T_{\tilde X} \oplus {f'}^* T_{\tilde X}
      \stackrel {\varphi}{\to} f_Q^* T_{\tilde X}
      \to 0 \\
    0 \to {f'}^* T_{\tilde X} (-Q)
      \to {f'}^* T_{\tilde X}
      \stackrel {\psi}{\to} f_Q^* T_{\tilde X}
      \to 0
  \end {gather*}
  where $ f_0 $, $ f' $, and $ f_Q $ denote the restrictions of $f$ to $ C_0 $,
  $ C' $, and $Q$, respectively.

  From these sequences we deduce that
  \begin {gather*}
    \dim \coker H^0 (\varphi)
      = h^1 (C, f^* T_{\bX}) - h^1 (C_0, f_0^* T_{\bX})
         - h^1 (C', {f'}^* T_{\bX}) \\
    \dim \coker H^0 (\psi)
      = h^1 (C', {f'}^* T_{\bX} (-Q)) - h^1 (C', {f'}^* T_{\bX}).
  \end {gather*}
  Since we certainly have $ \dim \coker H^0 (\varphi) \le \dim \coker H^0
  (\psi) $, we can combine these equations into the single inequality
    \[ h^1 (C, f^* T_{\bX})
         \le h^1 (C_0, f_0^* T_{\bX})
           + h^1 (C', {f'}^* T_{\bX} (-Q)). \]
  Now, by the induction hypothesis on $ f_0 $, we have $ h^1 (C_0, f_0^*
  T_{\bX}) \le \eta (C_0,f_0) $ with strict inequality holding if
  $ \eta (C_0,f_0) > 0 $. On the other hand, we get $ h^1 (C', {f'}^*
  T_{\bX} (-Q)) \le \eta (C',f') $ by lemma \ref {coh-irr} for $ \varepsilon=
  1 $. As $ \eta (C,f) = \eta (C_0,f_0) + \eta(C',f') $, the proposition
  follows by induction.
\end {proof}

We now come to the central proposition already alluded to in section \ref
{blowup_sign}: given a part $ \Mm {\bX,\tau} $ of the moduli space $ \Mb [n]%
{\bX,\beta} $ corresponding to the topology $ \tau $ (see section \ref
{blowup_prelim}), we consider the map
  \[ \phi : \Mm {\bX,\tau} \hookrightarrow \Mb [n]{\bX,\beta}
       \to \Mb [n]{X,d(\beta)} \]
given by mapping $ (C,x_1,\dots,x_n,f) $ to $ (C,x_1,\dots,x_n,p \circ f) $
and stabilizing if necessary ($ \phi $ exists by the functoriality of the
moduli spaces of stable maps, see \cite {BM} remark after theorem 3.14). We
show that, although $ \Mm {\bX,\tau} $ may have too big dimension, the image
$ \phi (\Mm {\bX,\tau}) $ has not. Part (ii) of the proposition, which is of
similar type, will be needed later in section \ref {blowup_tangency}.

\begin {proposition} \label {dim-im}
  Let $ \bX=\bP^r(1) $ and $ \beta \in A_1(\bX) $ with $ d(\beta) > 0 $.
  Let $ \phi: \Mb [n]{\bX,\beta} \to \Mb [n]{X,d(\beta)} $ be the morphism as
  above, and let $ \tau $ be a topology of stable maps of homology class
  $ \beta $ (so that $ \Mm {\bX,\tau} \subset \Mb [n]{\bX,\beta} $). Then we
  have
  \begin {enumerate}
  \item $ \dim \phi(\Mm {\bX,\tau}) \le \vdim \Mb [n]{\bX,\beta} $. Moreover,
    strict inequality holds if and only if $\tau$ is a topology corresponding
    to reducible curves.
  \item At least one of the following holds:
    \begin {itemize}
    \item [(a)] $ \dim \phi(\Mm {\bX,\tau}) \le \vdim \Mb [n]{\bX,\beta} -r $,
    \item [(b)] $ \dim \Mm {\bX,\tau} \le \vdim \Mb [n]{\bX,\beta} -2 $,
    \item [(c)] $ \dim \Mm {\bX,\tau} \le \vdim \Mb [n]{\bX,\beta} -1 $ and
      $ \eta(\tau) = 0 $,
    \item [(d)] $ \dim \Mm {\bX,\tau} = \vdim \Mb [n]{\bX,\beta} $ and $\tau$
      is the topology corresponding to irreducible curves,
    \item [(e)] $ \dim \Mm {\bX,\tau} = \vdim \Mb [n]{\bX,\beta} -1 $ and
      $\tau$ is a topology corresponding to reducible curves having exactly
      two irreducible components, one with homology class $ \beta-E' $ and the
      other with homology class $ E' $.
    \end {itemize}
  \end {enumerate}
\end {proposition}

\begin {proof}
  We start by defining some numerical invariants of the topology $ \tau $ that
  will be needed in the proof.
  \begin {itemize}
  \item Let $ \mdf {S} $ be the number of nodes of a curve with topology
    $\tau$. We divide this number into $ S = S_{EE} + S_{XX} + S_{XE} $,
    where $ \mdf {S_{EE}} $ (resp.\ $ \mdf {S_{XX}} $, $ \mdf {S_{XE}} $)
    denotes the number of nodes joining two exceptional components of $C$
    (resp.\ two non-exceptional components, or one exceptional with one
    non-exceptional component). Here and in the following we call an
    irreducible component of $C$ exceptional if it is mapped by $f$ into the
    exceptional divisor and it is not contracted by $f$, and non-exceptional
    otherwise.
  \item Let $ \mdf{P} $ be the (minimal) number of additional marked points
    which are necessary to stabilize $C$. We divide the number $P$ into $ P =
    P_E + P_X $, where $ \mdf {P_E} $ (resp.\ $ \mdf {P_X} $) is the number
    of marked points that have to be added on exceptional components (resp.\
    non-exceptional components) of $C$ to stabilize $C$.
  \end {itemize}
  Now we give an estimate for the dimension of $ \Mm {\bX,\tau} $. The tangent
  space $ T_{\Mm {\bX,\tau},\cc} $ at a point $ \cc = (C,x_1,\dots,x_n,f) \in
  \Mm {\bX,\tau} $ is given by the hypercohomology group (see \cite {K} section
  1.3.2)
    \[ T_{\Mm {\bX,\tau},\cc} = \HH^1 (T'_C \to f^* T_{\bX}) \]
  where $ T'_C = T_C (- x_1 - \dots - x_n) $ and where we put the sheaves
  $ T'_C $ and $ f^* T_{\bX} $ in degrees 0 and 1, respectively. This
  means that there is an exact sequence
  \begin {gather*}
    0 \to H^0 (C,T'_C)
      \to H^0 (C,f^* T_{\bX})
      \to T_{\Mm {\bX,\tau},\cc}
      \to H^1 (C,T'_C)
    \tag {1}
  \end {gather*}
  (note that the first map is injective because $ \cc $ is a stable map). By
  lemma \ref {coh-red}, we have
  \begin {gather*}
    \dim H^0 (C, f^* T_{\bX}) \le \chi (C, f^* T_{\bX}) + \eta (C,f).
       \tag {2}
  \end {gather*}
  Moreover, by Riemann-Roch we get $ \chi (C,T'_C) = S + 3 - n $. It follows
  that
  \begin {align*}
    \dim T_{\Mm {\bX,\tau},\cc} &\le \chi (C, f^* T_{\bX}) + \eta (C,f)
      + n - S - 3 \\
      &= \vdim \Mb [n]{\bX,\beta} + \eta(C,f) - S,
  \end {align*}
  and therefore
    \[ \dim \Mm {\bX,\tau} \le \vdim \Mb [n]{\bX,\beta} + \eta(\tau) - S. \]
  If $ \eta(\tau) - S < 0 $, then statement (i) is obviously satisfied.
  Moreover, if $ \eta(\tau)=0 $ then we also have (ii)-(c), and if $ \eta
  (\tau)>0 $ then we have strict inequality also in (2) and therefore (ii)-(b).
  Therefore we can assume from now on that $ \underline {\eta(\tau) - S \ge 0}
  $. If $ \eta(\tau) = 0 $, then we must also have $ S=0 $, which means that
  the curve is irreducible. But then (i) and (ii)-(d) are satisfied. So we can
  also assume in the sequel that $ \underline {\eta(\tau)>0} $. It follows
  then from lemma \ref {coh-red} that we have strict inequality in (2), hence
  \begin {gather*}
    \dim T_{\Mm {\bX,\tau},\cc} \le 
         \vdim \Mb [n]{\bX,\beta} + \eta(C,f) - S - 1.
    \tag {3}
  \end {gather*}
  We now give an estimate of the dimension of the image $ \phi (\Mm {\bX,\tau})
  $. As we always work over the ground field $ \CC $, we can do this on the
  level of tangent spaces, i.e.\ we have
    \[ \dim \phi (\Mm {\bX,\tau}) \le
         \max_{\cc \in \Mm {\bX,\tau}} \dim (d\phi)(T_{\Mm {\bX,\tau},\cc}). \]
  Hence our goal is to find as many vectors in $ \ker d\phi $ as possible. We
  do this by finding elements in the kernel of the composite map (see (1))
    \[ H^0 (C,f^* T_{\bX}) / H^0 (C,T'_C) \hookrightarrow T_{\Mm {\bX,\tau},
         \cc} \to T_{\Mb [n]{X,d(\beta)},\phi(\cc)}. \]
  Let $ C_0 $ be a maximal connected subscheme of $C$ consisting only of
  exceptional components of $C$. Let $ f_0 $ be the restriction of $f$ to
  $ C_0 $ and let $ Q_1,\dots,Q_a $ be the nodes of $C$ which join $C_0$ with
  the rest of $C$ (they are of type $ S_{XE} $). Now every section of $ f_0^*
  T_E (-Q_1- \cdots - Q_a) $ can be extended by zero to a section of $ f^*
  T_{\bX} $ which is mapped to zero by $ d\phi $ since these deformations
  of the map take place entirely within the exceptional divisor. As $ E \cong
  \PP^{r-1} $ is a convex variety, we have
    \[ h^0 (C_0, f_0^* T_E) = \chi (C_0, f_0^* T_E)
       = r-1 + r \, \eta (C_0,f_0) \]
  and therefore we can estimate the dimension of the space of deformations
  that we have just found:
    \[ h^0 (C_0, f_0^* T_E (-Q_1- \cdots -Q_a))
       \ge r-1 + r \, \eta (C_0,f_0) - (r-1) \, a. \]
  (The right hand side of this inequality may well be negative, but
  nevertheless the statement is correct also in this case, of course.)

  We will now add up these numbers for all possible $ C_0 $, say there are
  $ \mdf {B} $ of them. The sum of the $ \eta (C_0,f_0) $ will then give $
  \eta (C,f) $, and the sum of the $a$ will give $ S_{XE} $. Note that there
  is a $ P_E $-dimensional space of infinitesimal automorphisms of $C$, i.e.\
  a subspace of $ H^0 (C,T'_C) $, included in the deformations that we have
  just found, and that these are exactly the trivial elements in the kernel of
  $ d\phi $. Therefore we have
  \begin {align*}
    \dim \ker d\phi
      &\ge B \, (r-1) + r \, \eta (C,f) - (r-1) S_{XE} - P_E \\
      &=(r-2) (\underbrace {B\mathstrut}_{\ge 1}
        + \underbrace {\eta (C,f) - S_{XE}}_{\ge 0})
        + B + 2\eta (C,f) - S_{XE} - P_E \\
      & \qquad \quad \mbox {($ B \ge 1 $ since $ \eta(C,f)>0 $} \\
      & \qquad \quad \mbox {~~and $ \eta (C,f) - S_{XE} \ge 0 $ since $ \eta
        (C,f) - S \ge 0 $)} \\
      &\ge (r-2) + B + 2\eta (C,f) - S_{XE} - P_E.
  \end {align*}
  Combining this with (3), we get the estimate
  \begin {align*}
    \dim \phi(\Mm {\bX,\tau}) &\le \dim T_{\Mm {\bX,\tau},\cc} - \dim \ker
      d\phi \\
    &\le \vdim \Mb [n]{\bX,\beta} - r + 1 - (S_{XX}+S_{EE}+B+\eta(\tau)-P_E).
  \end {align*}
  To prove the proposition, it remains to look at the term in brackets. First
  we will show that
  \begin {align*}
    P_E \le S_{XX}+S_{EE}+B+\eta(\tau). \tag {4}
  \end {align*}
  Look at $ P_E $, i.e.\ the exceptional components of $C$ where marked points
  have to be added to stabilize $C$. We have to distinguish three cases:
  \begin {itemize}
    \item [(A)] Components on which two points have to be added, and whose
      (only) node is of type $ S_{EE} $: those give a contribution of 2 to
      $ P_E $, but they also give at least 1 to $ \eta (\tau) $ and to $
      S_{EE} $ (and every node of type $ S_{EE} $ belongs to at most one such
      component).
    \item [(B)] Components on which two points have to be added, and whose
      (only) node is of type $ S_{XE} $: those give a contribution of 2 to
      $ P_E $, but they also give at least 1 to $ \eta (\tau) $ and to $B$
      (since such a component alone is one of the $ C_0 $ considered above).
    \item [(C)] Components on which only one point has to be added: those
      give a contribution of 1 to $ P_E $, but they also give at least 1 to
      $ \eta (\tau) $.
  \end {itemize}
  This shows (4), finishing the proof of (i). As for (ii), (a) is satisfied if
  we have strict inequality in (4), so we assume from now on that this is not
  the case and determine necessary conditions for equality by looking at the
  proof of (4) above. First of all, we see that every maximal connected
  subscheme of $C$ consisting only of exceptional components contributes 1 to
  $B$, but this gets accounted for only in case (B) above, so if we want to
  have equality, every such maximal connected subscheme must actually be an
  irreducible component of type (B), which in addition gives a contribution of
  \emph {exactly} 2 to $ P_E $ and \emph {exactly} 1 to $ \eta (\tau) $. So all
  exceptional components of the curve must actually be lines with no marked
  points, connected at exactly one point to a non-exceptional component of
  the curve. Moreover, for equality we must also have $ S_{XX} = 0 $, since
  these nodes have not been considered above at all.

  Hence, in summary, we must have one non-exceptional irreducible component
  $ C_0 $ of homology class $ \beta - q\,E' $, and $q$ exceptional components
  of homology class $ E' $ with no marked points, each connected at exactly
  one point to $ C_0 $. But it is easy to compute the dimension of $ \phi
  (\Mm {X,\tau}) $ for these topologies: the map $ \phi $ simply forgets the
  $q$ exceptional components, so
  \begin {align*}
    \dim \phi (\Mm {\bX,\tau}) &= \dim \Mm [n]{\bX,\beta-q\,E'} \\
      &= \vdim \Mb [n]{\bX,\beta-q\,E'} \qquad \qquad \mbox {by (i)} \\
      &= \vdim \Mb [n]{\bX,\beta} - q(r-1).
  \end {align*}
  Hence we see that (ii)-(a) is satisfied for $ q>1 $ and (ii)-(e) for
  $ q=1 $.

  This completes the proof.
\end {proof}

We now combine our results to prove the enumerative significance of the
\gwinvs of $ \bP^r(1) $. Some examples of these numbers can be found in
\ref {p2-1} and \ref {p3-1}.

\begin {theorem} \label {enum-1}
  Let $ \bX = \bP^r(1) $, $ \beta = d\,H'+e\,E' \in A_1(\bX) $ an effective
  homology class with $ d>0 $ and $ e \le 0 $, and $ \calT = \gamma_1 \seq
  \dots \seq \gamma_n $ a collection of non-exceptional effective classes such
  that $ \sum_i \codim \gamma_i = \vdim {\Mb [n]{\bX,\beta}} $. Then $ \gw
  {\beta}{\calT} $ is enumerative.
\end {theorem}

\begin {proof}
  The proof goes along the same lines as that of lemma \ref {bXfromX}. For
  irreducible stable maps $ (C,x_1,\dots,x_n,f) $ we have $ h^1 (C,f^*
  T_{\bX}) = 0 $ by lemma \ref {coh-irr} (i). Therefore, if $ Z \subset \Mb
  [n]{\bX,\beta} $ denotes the closure of $ \Mm [n]{\bX,\beta} $, then lemma
  \ref {vfc} tells us that
    \[ \virt {\Mb [n]{\bX,\beta}} = [Z] + \alpha \]
  where $ \alpha $ is a cycle of dimension $ \vdim \Mb [n]{\bX,\beta} $
  supported on $ \Mb [n]{\bX,\beta} \backslash \Mm [n]{\bX,\beta} $. But if
  $ \phi: \Mb [n]{\bX,d\,H'+e\,E'} \to \Mb [n]{X,d\,H'} $ denotes the morphism
  induced by the map $ p: \bX \to X $, we must have $ \phi_* \alpha = 0 $ by
  proposition \ref {dim-im} (i). So, considering the commutative diagram
  \xydiag {
    \Mb [n]{\bX,\beta} \ar[r]^{\phi} \ar[d]_{ev_i}
      & \Mb [n]{X,d\,H'} \ar[d]_{ev_i} \\
    \bX \ar[r]^p & X
  }
  for $ 1 \le i \le n $, it follows by the projection formula that
  \begin {align*}
    \gw [\bX]{\beta}{\calT}
      &= (\prod_i ev_i^* p^* \gamma_i) \cdot \virt {\Mb [n]{\bX,\beta}} \\
      &= (\prod_i ev_i^* \gamma_i) \cdot \phi_* \virt {\Mb [n]{\bX,\beta}} \\
      &= (\prod_i ev_i^* \gamma_i) \cdot \phi_* [Z]. \\
      &= (\prod_i ev_i^* p^* \gamma_i) \cdot [Z].
  \end {align*}
  Hence we are evaluating an intersection product on the stack $Z$.

  Unless $ d+e=0 $ and $ d\ge 2 $, the theorem now follows from the
  Bertini lemma \ref {bertini} (iii) in combination with lemma \ref {noauto}
  saying that the generic element of $Z$ has no automorphisms and corresponds
  to a generically injective stable map. However, if $ d+e=0 $ and $ d \ge 2
  $, then the image of every stable map in $ \Mm [n]{\bX,d\,H'-d\,E'} $ is a
  line through the blown-up point. These curves can obviously only satisfy as
  many incidence conditions as the curves in $ \Mm [n]{\bX,H'-E'} $. But $
  \vdim \Mb [n]{\bX,d\,H'-d\,E'} > \vdim \Mb [n]{\bX,H'-E'} $, so the \gwinv
  will be zero, which is also the enumeratively correct number.
\end {proof}


\section {Enumerative significance --- the case $ \bP^3(4) $} \label
  {blowup_four}

In this section, we discuss the enumerative significance of the \gwinvs on
$ \bX=\bP^3(4) $. First we fix some notation. As the four points to blow up
on $ X=\PP^3 $ we choose $ P_1 = (1:0:0:0) $, $ P_2 = (0:1:0:0) $, $ P_3 =
(0:0:1:0) $, and $ P_4 = (0:0:0:1) $. For $ 1 \le i<j \le 4 $, we denote by
$ L_{ij} \subset \bX $ the strict transform of the line $ \overline {P_i P_j}
$. The $ L_{ij} $ are disjoint from each other, and we set $ \calL =
\bigcup_{i<j} L_{ij} $. For $ 1 \le i \le 4 $, we let $ H_i $ be the strict
transform of the hyperplane in $X$ spanned by the three points $ P_j $ with
$ j \neq i $, and we set $ \calH = \bigcup_i H_i $. As usual, $ E_i $
denotes the exceptional divisor over $ P_i $. We set $ \calE = \bigcup_i E_i $.

Let $ \beta \in A_1(\bX) $ be an effective homology class with $ d(\beta) >
0 $. The first thing to do is to look at \emph {irreducible} curves of
homology class $ \beta $ and to see whether their moduli space $ \Mm [0]{
\bX,\beta} $ is smooth and of the expected dimension, which in this case is
  \[ \vdim \Mb [0]{\bX,\beta} = 4d(\beta)+2e(\beta). \]
In the case of one \blowup in section \ref {blowup_one}, this followed easily
from lemma \ref {coh-irr} (i) since there we always have $ d(\beta)+e(\beta)
\ge 0 $. However, for multiple \blowups, this is not necessarily the case. Our
way to solve this problem is to use a certain Cremona map to transform curves
with $ d(\beta)+e(\beta) \le 0 $ into others with $ d(\beta)+e(\beta) \ge 0 $,
so that lemma \ref {coh-irr} can be applied again. Before we can describe this
map, we need a definition.

\begin {definition}
  Let $ (C,f) \in \Mm [0]{\bP^3(4),\beta} $ be an irreducible stable
  map with $ f(C) \not\subset \calL $. Then we set $ \lambda_{ij} (C,f) $
  to be the ``multiplicity of $f$ along $ L_{ij} $'', defined as follows:
  if $ \varphi_1: \tilde Y \to \bP^3 (4) $ is the \blowup of $ \bP^3 (4) $
  along $ \calL $ with exceptional divisors $ F_{ij} $ over $ L_{ij} $, then
  there is a well-defined map $ \varphi_1^{-1} \circ f: C \to \tilde Y $, and
  we define
    \[ \mdf {\lambda_{ij} (C,f)} := F_{ij} \cdot (\varphi_1^{-1} \circ f)_*
         [C] \ge 0. \]
  Finally, we define $ \mdf {\vec\lambda (C,f)} $ to be the vector consisting
  of all $ \lambda_{ij} (C,f) $, and set
    \[ \mdf {\lambda (C,f)} = \sum_{i<j} \lambda_{ij}(C,f). \]
\end {definition}

We can now describe the Cremona map announced above.

\begin {lemma} \label {cremona}
  There exists a birational map $ \varphi: \bP^3(4) \dashrightarrow \bP^3(4) $
  which is an isomorphism outside $ \calL $ with the following property:

  If $ (C,f) \in \Mm [0]{\bP^3(4),\beta} $ is an irreducible stable
  map such that $ f(C) \not\subset \calL $, so that the transformed stable map
  $ (C,\varphi \circ f) \in \Mm [0]{\bP^3(4),\beta'} $ exists, then the
  homology class $ \beta' $ of the transformed stable map satisfies
  \begin {align*}
    d(\beta') &= 3d(\beta)+2e(\beta)-\lambda(C,f), \\
    e(\beta') &= -4d(\beta)-3e(\beta)+2\lambda(C,f).
  \end {align*}
  Hence, in particular, we have
  \begin {itemize}
  \item $ 4d(\beta')+2e(\beta') = 4d(\beta)+2e(\beta) $,
  \item if $ d(\beta)+e(\beta) \le 0 $, then $ d(\beta')+e(\beta') \ge 0 $.
  \end {itemize}
\end {lemma}

\begin {proof}
  The birational map $ \varphi: \bP^3(4) \dashrightarrow \bP^3(4) $ we want to
  consider is most easily described in the language of toric geometry (see
  e.g.\ \cite {F}). Let $ \Delta' $ in $ \RR^3 $ be the complete simplicial
  fan with one-dimensional cones $ \{ \langle v_i \rangle \;|\; 1 \le i \le
  4 \} $, where 
    \[ v_1 = (1,0,0),\; v_2 = (0,1,0),\; v_3 = (0,0,1),\; v_4 = (-1,-1,-1), \]
  corresponding to the toric variety $ X_{\Delta'} = \PP^3 $. Let $ \Delta $
  be the \blowup of $ \Delta' $ at the four torus-invariant points as described
  in \cite {F} section 2.4, so that the toric variety $ X_\Delta $
  associated to $ \Delta $ is $ \bP^3(4) $. The fan $ \Delta $ can be
  described explicitly as follows: it is the complete fan with one-dimensional
  cones
    \[ \{ \pm \langle v_i \rangle \;|\; 1 \le i \le 4 \} \]
  and two-dimensional cones
    \[ \{ \langle v_i,-v_j \rangle \;|\; 1 \le i,j \le 4 ;\, i \neq j \}
         \cup \{ \langle v_i, v_j \rangle \;;\; 1 \le i<j \le 4 \}. \]
  The Picard group of $ X_\Delta $ is generated by the divisors corresponding
  to the one-dimen\-sional cones, we will denote the divisor corresponding to
  the cone $ \langle v_i \rangle $ by $ H_i $ and the divisor corresponding to
  the cone $ - \langle v_i \rangle $ by $ E_i $. This coincides with the
  definition of $ H_i $ and $ E_i $ given above, and these divisors satisfy
  the three relations
  \begin {align*}
    H :&= H_1 + E_2 + E_3 + E_4 \\
       &= H_2 + E_1 + E_3 + E_4 \\
       &= H_3 + E_1 + E_2 + E_4 \\
       &= H_4 + E_1 + E_2 + E_3 \tag {1}
  \end {align*}
  where $H$ denotes the pullback of the hyperplane class under the map
  $ p: \bP^3(4) \to \PP^3 $.

  Now denote by $ -\Delta $ the fan obtained by mirroring $ \Delta $ at the
  origin in $ \RR^3 $. Then, of course, we also have $ X_{-\Delta} \isom
  \bP^3(4) $. The map $ \varphi $ we want to consider is now the obvious
  rational map $ \varphi: X_\Delta \dashrightarrow X_{-\Delta} $ which is the
  identity on the torus $ (\CC^*)^3 $ contained in both $ X_\Delta $ and $
  X_{-\Delta} $. Note that the one-dimensional cones of $ \Delta $ and $
  -\Delta $ are the same, so that $ \varphi $ is an isomorphism away from a
  subvariety of $ \bP^3(4) $ of codimension 2.

  In more geometric terms, we can describe $ \varphi $ as the so-called
  ``flip'' of the 6 lines $ \calL $, i.e.\ one blows up these lines
  (that have normal bundle $ \OO (-1) \oplus \OO(-1) $ in $ \bP^3 (4) $) to
  get a variety $ \tilde Y $ with the 6 exceptional divisors $ \hat F_{ij}
  \isom \PP^1 \times \PP^1 $ corresponding to $ L_{ij} $, and then blows down
  the $ F_{ij} $ again with the roles of base and fibre reversed in $ \PP^1
  \times \PP^1 $. One can write these two steps as in the following diagram:
  \xydiag {
    & \tilde Y \ar[dl]_{\varphi_1} \ar[dr]^{\varphi_2} \\
    X_\Delta \isom \bP^3(4) \ar@{-->}[rr]^{\varphi} & &
      X_{-\Delta} \isom \bP^3(4).
  }
  The variety $ \tilde Y $ can be depicted as follows:
  \begin {center} \unitlength 1cm \begin {picture}(10,7)
    \put (0,0){\epsfig {file=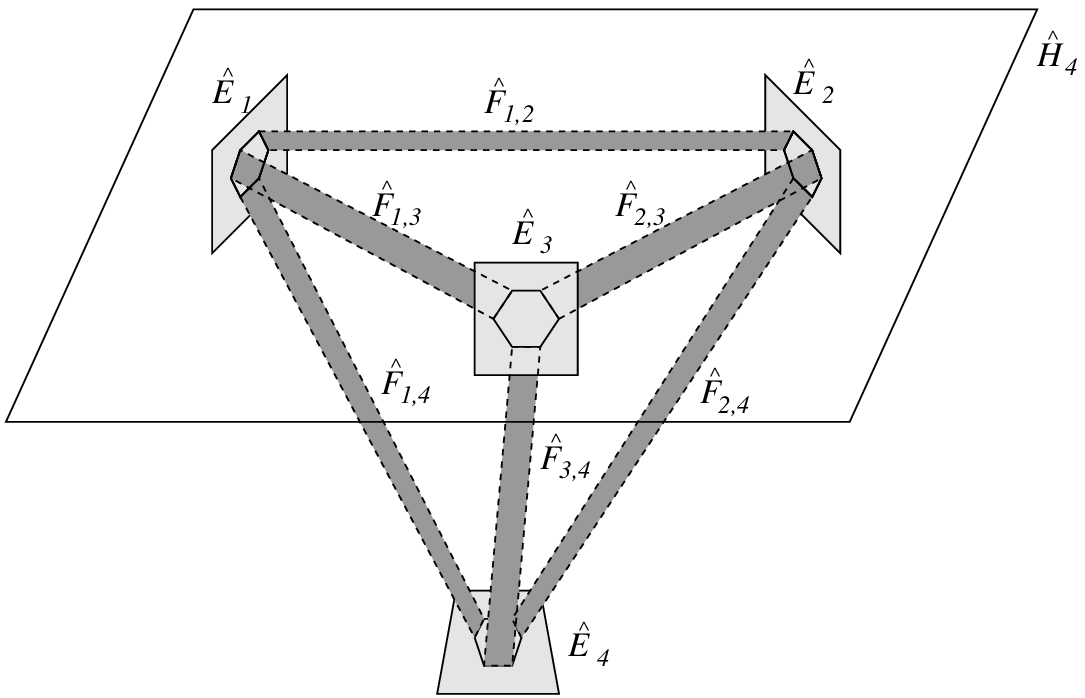,width=10cm}}
  \end {picture} \end {center}
  Here, we denoted the strict transforms of $ H_i $ and $ E_i $ under
  $ \varphi_1 $ by $ \hat H_i $ and $ \hat E_i $, respectively. These are
  all isomorphic to $ \bP^2(3) $. The divisors $ \hat H_1 $, $ \hat H_2 $,
  and $ \hat H_3 $ have not been drawn to keep the picture simple.

  We now look more closely at the divisors in $ \tilde Y $. Obviously, we have
  \begin {align*}
    \varphi_1^* H_1 &= \hat H_1 + \hat F_{23} + \hat F_{24} + \hat F_{34}, \\
    \varphi_1^* E_1 &= \hat E_1,
  \end {align*}
  and similarly for $ H_i $ and $ E_i $ with $ i=2,3,4 $. The Picard group of
  $ \tilde Y $ is the free abelian group generated by the $ \hat H_i $, $
  \hat E_i $, and $ \hat F_{ij} $, modulo the three relations induced by (1)
  \begin {align*}
    \hat H := \varphi_1^* H
       &= \hat H_1 + \hat E_2 + \hat E_3 + \hat E_4
            + \hat F_{23} + \hat F_{24} + \hat F_{34} \\
       &= \hat H_2 + \hat E_1 + \hat E_3 + \hat E_4
            + \hat F_{13} + \hat F_{14} + \hat F_{34} \\
       &= \hat H_3 + \hat E_1 + \hat E_2 + \hat E_4
            + \hat F_{12} + \hat F_{14} + \hat F_{24} \\
       &= \hat H_4 + \hat E_1 + \hat E_2 + \hat E_3
            + \hat F_{12} + \hat F_{13} + \hat F_{23}. \tag {2}
  \end {align*}
  If we now have a stable map in $ (C,f) $ in $ \tilde Y $, we also get stable
  maps $ (C_i,f_i) $ in $ \bP^3(4) $ by composing $f$ with $ \varphi_i $ for
  $ i=1,2 $. We will now compute the homology classes of these two stable maps.

  The homology class of $ (C_1,f_1) $ is $ \beta = d\,H'+\sum_i e_i\,E'_i $
  where
  \begin {align*}
    d &= H \cdot {\varphi_1}_* f_* [C] \\
      &= \hat H \cdot f_* [C] \\
      &= (\hat H_1 + \hat E_2 + \hat E_3 + \hat E_4
           + \hat F_{23} + \hat F_{24} + \hat F_{34}) \cdot f_* [C], \\
    e_i &= - E_i \cdot {\varphi_1}_* f_* [C] \\
        &= - \hat E_i \cdot f_* [C].
  \end {align*}
  The homology class of $ (C_2,f_2) $ is obtained by reversing the roles of
  $ \hat H_i $ and $ \hat E_i $ and substituting $ \hat F_{12} \leftrightarrow
  \hat F_{34} $, $ \hat F_{13} \leftrightarrow \hat F_{24} $, and $ \hat
  F_{14} \leftrightarrow \hat F_{23} $, so it is $ \beta' = d'\,H'+\sum_i
  e'_i\,E'_i $ where
  \begin {align*}
    d' &= (\hat E_1 + \hat H_2 + \hat H_3 + \hat H_4
           + \hat F_{14} + \hat F_{13} + \hat F_{12}) \cdot f_* [C] \\
       &= (3\hat H_1 - 2\hat E_1 + \hat E_2 + \hat E_3 + \hat E_4
           - \hat F_{12} - \hat F_{13} - \hat F_{14}
           + 2\hat F_{23} + 2\hat F_{24} + 2\hat F_{34}) \cdot f_* [C] \\
       &  \qquad \quad \mbox {(by substituting $ \hat H_2 $, $ \hat H_3 $, and
            $ \hat H_4 $ from (2))} \\
       &= 3d+2(e_1+e_2+e_3+e_4)-
            \underbrace {(\sum_{i<j} F_{ij}) \cdot f_* [C]}_{
            =\lambda(C_1,f_1)=\lambda(C_2,f_2)=:\lambda}, \\
    e'_1 &= - \hat H_1 \cdot f_* [C] \\
         &= - d - e_2 - e_3 - e_4 + (\hat F_{23} + \hat F_{24} + \hat F_{34})
              \cdot f_* [C],
  \end {align*}
  and similarly for $ e_2 $, $ e_3 $, and $ e_4 $. Defining $ e = \sum_i e_i $
  and $ e' = \sum_i e'_i $, we arrive at the equations
  \begin {align*}
    d' &= 3d+2e-\lambda, \\
    e' &= -4d-3e+2\lambda.
  \end {align*}
  In particular, we see that $ 4d'+2e' = 4d+2e $ and that, if $ d+e \le 0 $,
  then
    \[ d'+e' = -d-e+\lambda \ge \lambda \ge 0. \]
\end {proof}

We now use this map to prove some properties of irreducible stable maps in
$ \bX=\bP^3(4) $. As already mentioned in section \ref {blowup_sign}, apart
from the case where $ \Mm [n]{\bX,\beta} $ is smooth of the expected dimension
(case (iii) below), we have to consider the cases where the curves are
multiple coverings of one of the $ L_{ij} $ (case (i)) and where they
are contained in one of the $ H_i $ (such that they cannot satisfy any
incidence conditions with generic points in $ \bX $, see case (ii)). One of
the most important statements of the next lemma is the final conclusion that,
although the dimension of the moduli space may be too big, the curves can
never satisfy more incidence conditions (with points) as one would expect
from the virtual dimension of the moduli space.

\begin {lemma} \label {p3-prep}
  Let $ \beta \in A_1 (\bX) $ be a homology class such that $ \Mm [0]{\bX,
  \beta} \neq \emptyset $. Set
    \[ n := \frac 1 2 \vdim \Mb [0]{\bX,\beta} = 2 d(\beta)+e(\beta). \]
  Then at least one of the following statements holds:
  \begin {enumerate}
  \item $ n=0 $ and $ \beta = d\,H'-d\,E'_i-d\,E'_j $ for some $ d>0 $,
    $ 1 \le i < j \le 4 $. All curves in $ \Mm [0]{\bX,\beta} $ are contained
    in $ L_{ij} $.
  \item $ n>0 $, and for generic points $ Q_1,\dots,Q_n \in \bX $, we have
      \[ ev_1^{-1} (Q_1) \cap \dots \cap ev_n^{-1} (Q_n) = \emptyset \]
    in $ \Mm [n]{\bX,\beta} $.
  \item $ n>0 $, $ \dim \Mm [0]{\bX,\beta} = \vdim \Mb [0]{\bX,\beta} $,
    and for a generic element $ \cc = (C,f) \in \Mm [0]{\bX,\beta} $, $f$ is
    generically injective, $ \cc $ has no automorphisms, and $ f(C) $
    intersects neither $ \calL $ (which is a disjoint union of 6 smooth
    rational curves) nor $ \calH \cap \calE $ (which is a union of 12 smooth
    rational curves).
  \end {enumerate}
  In particular, it is impossible that $ n<0 $, and in any case we have
    \[ ev_1^{-1} (Q_1) \cap \dots \cap ev_{n'}^{-1} (Q_{n'}) = \emptyset \]
  in $ \Mm [n']{\bX,\beta} $ for generic points $ Q_1,\dots,Q_{n'} \in \bX $
  if $ n' > n $.
\end {lemma}

\begin {proof}
  Let $ (C,f) \in \Mm [0]{\bX,\beta} $ be a stable map, $ d=d(\beta) $, $
  e_i = e_i(\beta) $, $ e=\sum_i e_i $, and assume that $ \beta \neq 0 $
  (since otherwise $ \Mm[0]{\bX,\beta} = \emptyset $).

  If $ d=0 $, then $ n = e(\beta) > 0 $ and $ f(C) $ is contained in an
  exceptional divisor. Then it is clear that for a generic point in $ \bX $,
  no curve in $ \Mm [0]{\bX,\beta} $ meets this point. Therefore, (ii) is
  satisfied.

  Now assume $ d>0 $, then we must have $ e_i \le 0 $ for all $i$. The curve
  $ f(C) $ cannot be contained at the same time in three of the $ H_i $, since
  their intersection is empty. This means that there are at least two of the
  $ H_i $, say $ H_1 $ and $ H_2 $, in which $ f(C) $ is not contained. It
  follows that
    \[ d+e_2+e_3+e_4 = \deg f^* H_1 \ge 0 \quad \mbox {and} \quad
       d+e_1+e_3+e_4 = \deg f^* H_2 \ge 0. \]
  Since $ e_4 \le 0 $ and $ e_3 \le 0 $, this also means that $ d+e_2+e_3 \ge
  0 $ and $ d+e_1+e_4 \ge 0 $, and therefore $ n = 2d+e \ge 0 $: the virtual
  dimension of the moduli space cannot be negative. Moreover, if $ n=0 $ then
  we must have equality everywhere, which means
    \[ e_1=-d, \; e_2=-d, \; e_3=0, \; e_4=0. \]
  Hence we are in case (i), and it is clear that all these curves are $d$-%
  fold coverings of $ L_{12} $.

  It remains to consider the case when $ n>0 $. We distinguish four cases.

  \underline {Case 1: $ \beta = d\,H'-d\,E'_i $ for $ d>1 $ and some $ 1
  \le i \le 4 $.} Then the curves in $ \Mm [0]{\bX,\beta} $ must obviously
  be $d$-fold coverings of a line through the exceptional divisor $ E_i $.
  Those cannot pass through two generic points, however $ n=2d-d=d \ge 2 $,
  hence (ii) is satisfied.

  We assume therefore from now on that $ \beta $ is not of this form.

  \underline {Case 2: $ d+e \ge 0 $.} We show that (iii) is satisfied.
  \begin {itemize}
  \item $ \dim \Mm [0]{\bX,\beta} = \vdim \Mb [0]{\bX,\beta} $: This follows
    because $ h^1 (C,f^* T_{\bX}) = 0 $ by lemma \ref {coh-irr} (i).
  \item the generic element of $ \Mm [0]{\bX,\beta} $ has no automorphisms and
    corresponds to a generically injective map: This follows from lemma \ref
    {noauto}.
  \item the generic element of $ \Mm [0]{\bX,\beta} $ does not intersect
    $ \calL $ and $ \calH \cap \calE $: Let $L$ be one of the 18 smooth
    rational curves in $ \calL \cup (\calH \cap \calE) $, we will show that
    the generic element of $ \Mm [0]{\bX,\beta} $ does not intersect $L$.
    Assume that $ (C,f) $ is a stable map in $ \bX $ such that there is a
    point $ x \in C $ with $ f(x) = Q \in L $. Consider $ \cc = (C,x,f) $ as
    an element of $ M = \Mm [1]{\bX,\beta} $. The tangent space to $M$ at the
    point $ \cc $ is (see \cite {K} section 1.3.2)
      \[ T_{M,\cc} = H^0 (C,f^* T_{\bX}) / H^0 (C, T_C(-x)). \]
    If $ Z \subset M $ denotes the substack of those stable maps with
    $ f(x) \in L $, then the tangent space to $Z$ at $\cc$ is
      \[ T_{Z,\cc} = \{ s \in T_{M,\cc} \;;\; s(x) \in f^* T_{L,Q} \}. \]
    However, by lemma \ref {coh-irr} (i) for $ \varepsilon=1 $ we see that
      \[ h^0 (C,f^* T_{\bX} (-x)) = h^0 (C,f^* T_{\bX}) - 3, \]
    i.e.\ that the map $ H^0 (C,f^* T_{\bX}) \to f^* T_{\bX,Q}, \; s \mapsto
    s(x) $ is surjective. Therefore the tangent space to $Z$ at $\cc$ has
    smaller dimension than that to $M$. Since $M$ is smooth at $\cc$, it
    follows that $Z$ has smaller dimension than $M$ at $\cc$, proving the
    statement that the generic element of $ \Mm [0]{\bX,\beta} $ does not
    intersect $L$.
  \end {itemize}

  \underline {Case 3: $ d+e < 0 $ and $ e_i = 0 $ for some $i$.}
  Without loss of generality assume that $ e_4 = 0 $. Since then $ 0 > d+e =
  \deg f^* (H-E_1-E_2-E_3) = \deg f^* H_4 $, we conclude that $ f(C) $ must be
  contained in $ H_4 $. Hence (ii) is satisfied.

  \underline {Case 4: $ d+e < 0 $ and all $ e_i \neq 0 $.} We show that (iii)
  is satisfied using the Cremona map of lemma \ref {cremona}. We use in the
  following proof the notations of this lemma. Certainly no curve in $ \Mm
  [0]{\bX,\beta} $ is contained in $ \calL $. So if we decompose $ \Mm[0]{
  \bX,\beta} $ into parts $ M_{\vec \lambda} $ according to the value of $
  \vec \lambda(C) $ then $ \varphi $ gives injective morphisms from $ M_{\vec
  \lambda} $ to $ \Mm [0]{\bX,\beta_{\vec\lambda}} $ with $ \beta_{\vec\lambda
  } $ calculated in the proof of lemma \ref {cremona}. In particular we have
  $ d(\beta_{\vec\lambda}) + e(\beta_{\vec\lambda}) \ge 0 $, so that we can
  apply the results of case 2 to $ \Mm [0]{\bX,\beta_{\vec \lambda}} $. We
  therefore have
  \begin {align*}
    \dim M_{\vec\lambda} &\le \dim \Mm [0]{\bX,\beta_{\vec\lambda}}
           \tag {1} \\
      &= \vdim \Mb [0]{\bX,\beta_{\vec\lambda}} \quad
           \mbox {by case 2} \\
      &= 4d(\beta_{\vec\lambda}) + 2e(\beta_{\vec\lambda}) \\
      &= 4d(\beta) + 2e(\beta) \quad \mbox {by lemma \ref {cremona}} \\
      &= \vdim \Mb [0]{\bX,\beta}.
  \end {align*}
  If $ \vec\lambda \neq 0 $, i.e.\ if all curves in $ M_{\vec\lambda} $
  intersect $ \calL $, then the transformed curves in $ \Mm [0]{\bX,\beta_{
  \vec\lambda}} $ also have to intersect $ \calL $. But the generic curve in
  $ \Mm [0]{\bX,\beta_{\vec\lambda}} $ does not intersect $ \calL $ by the
  results of case 2, so it follows that we must have strict inequality in
  (1). Since the dimension of $ \Mb [0]{\bX,\beta} $ cannot be smaller than
  its virtual dimension, this means that $ M_{\vec\lambda} $ is nowhere dense
  in $ \Mm [0]{\bX,\beta} $ for $ \vec\lambda \neq \vec 0 $. In other words,
  $ M_{\vec 0} $ is dense in $ \Mm [0]{\bX,\beta} $, so it obviously suffices
  to prove (iii) for $ M_{\vec 0} $.

  But this is now easy: it follows from the above calculation that the
  dimension of $ M_{\vec 0} $ is equal to the virtual dimension of $ \Mb [0]{
  \bX,\beta} $. The other statements of (iii) about the generic curves in the
  moduli space are obviously preserved by the Cremona map $ \varphi $, so they
  follow from the fact that the space $ \Mm [0]{\bX,\beta_{\vec 0}} $ has these
  properties.

  This completes the proof that we always have one of the cases (i) to (iii).
  The statement that $ n \ge 0 $ has already been proven, and the fact that
    \[ ev_1^{-1} (Q_1) \cap \dots \cap ev_{n'}^{-1} (Q_{n'}) = \emptyset \]
  in $ \Mm [n']{\bX,\beta} $ for generic points $ Q_1,\dots,Q_{n'} \in \bX $
  if $ n' > n $ follows easily in all cases: for (i) because the image of
  all curves in the moduli space is contained in an $ L_{ij} $, for (ii) it
  is trivial, and for (iii) it follows from the Bertini lemma \ref {bertini}
  (ii).
\end {proof}

To prove enumerative significance for the \gwinvs on $ \bP^3(4) $, we now
finally have to consider reducible stable maps. Some numerical examples can
be found in \ref {p3-2}.

\begin {theorem} \label {enum-2}
  Let $ \bX = \bP^3(4) $ and $ \beta \in A_1(\bX) $ an effective homology
  class which is not of the form $ d\,H'-d\,E'_i-d\,E'_j $ for some $ d \ge 2
  $ and $ i \neq j $. Let $ \calT = \pt^{\seq n} $, where $ n=2d(\beta)+e(
  \beta) $. Then $ \gw {\beta}{\calT} $ is enumerative.
\end {theorem}

\begin {proof}
  Let $ Q_1,\dots,Q_n $ be generic points in $ \bX $. First we want to show
  that all points in the intersection
  \begin {gather*}
    I := ev_1^{-1} (Q_1) \cap \dots \cap ev_n^{-1} (Q_n) \tag {1}
  \end {gather*}
  on $ \Mb [n]{\bX,\beta} $ correspond to irreducible stable maps. To do this,
  we decompose the moduli space $ \Mb [n]{\bX,\beta} $ into the spaces $ M_\tau
  := \Mm {\bX,\tau} $ according to the topology of the curves and show that
  $ I \cap M_\tau $ is empty for each $ \tau $ corresponding to reducible
  curves.

  So assume that $ \tau $ is a topology corresponding to stable maps $ (C,f) $
  whose irreducible components \emph {that are not contracted by $f$} are
  $ C_1,\dots,C_a $. For $ 1 \le i \le a $, let $ \beta_i \neq 0 $ be the
  homology class of $f$ on $ C_i $ and let $ n_i $ be the number of markings
  on the component $ C_i $.

  By a \df {maximal contracted subscheme} we will mean a maximal connected
  subscheme of $C$ consisting only of components of $C$ that are contracted
  by $f$. A maximal contracted subscheme will be called \df {marked} if it
  contains at least one of the marked points. For each $ 1 \le i \le a $, we
  define $ \rho_i $ to be the number of marked maximal contracted subschemes
  of $C$ that have non-empty intersection with $ C_i $.

  We can assume that each maximal contracted subscheme hast at most one
  marked point, since otherwise the intersection (1) will certainly be
  empty. This means that each maximal contracted subscheme must have at
  least two points of intersection with the other components of the curve,
  since otherwise the prestable map $ (C,x_1,\dots,x_n,f) $ would not be
  stable. We conclude that each marked point that lies in a contracted
  component (there are $ (n-\sum_i n_i) $ of them) must be counted in at least
  two of the $ \rho_i $:
  \begin {gather*}
    \sum_i \rho_i \ge 2 (n-\sum_i n_i). \tag {2}
  \end {gather*}
  Now there is a morphism
  \begin {gather*}
    \Phi: M_\tau \to \Mm [n_1+\rho_1]{\bX,\beta_1} \times \dots \times
                  \Mm [n_a+\rho_a]{\bX,\beta_a} \tag {3}
  \end {gather*}
  mapping a stable map $ \cc $ to its non-contracted components, where on each
  such component we take as marked points the $ n_i $ marked points of $ \cc $
  lying on this component together with the intersection points of the
  component with the maximal contracted subschemes. We denote by $ \Phi_i:
  M_\tau \to \Mm [n_i+\rho_i]{\bX,\beta_i} $ the composition of $ \Phi $
  with the projections onto the factors of the right hand side of (3).

  We now consider again the intersection $I$ in (1) and show that $ \Phi
  (I \cap M_\tau) $ is empty for all topologies $ \tau $ but the trivial one,
  hence showing that $ I \cap M_\tau $ is empty. Note that in $ \Phi_i (I \cap
  M_\tau) $ the image point of each of the $ n_i+\rho_i $ marked points is
  fixed to be a certain $ Q_j $. But we have seen in lemma \ref {p3-prep}
  that, if $ \Phi_i (I \cap M_\tau) \subset \Mm [n_i+\rho_i]{\bX,\beta_i} $ is
  non-empty, this requires $ n_i+\rho_i $ to be at most $ 2d(\beta_i)+e(
  \beta_i) $. Therefore we get
  \begin {align*}
    n \le 2n-\sum_i n_i \stackrel {(2)}{\le} &\sum (n_i + \rho_i)
      \le \sum_i (2d(\beta_i)+e(\beta_i)) \\
    &= 2d(\beta)+e(\beta) = \frac 1 2 \vdim \Mb [0]{\bX,\beta} = n.
  \end {align*}
  Hence we must have equality everywhere, which means first of all that
  $ \sum_i n_i = n $ and therefore $ \rho_i = 0 $ for all $i$. Moreover,
  it follows that the number $ n_i $ of marked points with prescribed
  image in $ \Phi_i (I \cap M_\tau) $ is equal to $ 2d(\beta_i)+e(\beta_i) $
  for all $i$, showing that there can be no component of $C$ of type (ii)
  according to the classification of lemma \ref {p3-prep} (to be precise, that
  for all $i$, $C$ is mapped under $ \Phi_i $ to a moduli space which is not
  of type (ii)). If there are only components of type (i), then we have the
  case that $ \beta = d\,H-d\,E'_i-d\,E'_j $ for some $ d>2 $ and $ i \neq j $
  (note that there cannot be two components of type (i) with different $ (i,j)
  $ since the $ L_{ij} $ do not intersect). As we excluded this case in the
  theorem, we conclude that there must be at least one component of $C$ of
  type (iii). We are going to show that there is in fact only one component
  which must then necessarily be of type (iii).

  We first exclude components of type (i). Note that on each component $ C_i $
  of type (iii) we impose $ n_i $ generic point conditions. Since $ \dim
  \Mm [n_i]{\bX,\beta_i} = 3n_i $, this means by the Bertini lemma \ref
  {bertini} (ii) that $ \Phi_i (I \cap M_\tau) \subset \Mm [n_i]{\bX,\beta_i} $
  is zero-dimensional (if not empty). Moreover, if we let $ Z_i \subset
  \Mm [n_i]{\bX,\beta_i} $ be the substack of curves intersecting $
  \calL \cup (\calH \cap \calE) $, then $ \dim Z_i < 3n_i $ by lemma \ref
  {p3-prep}, and hence again by Bertini, $ \Phi_i (I \cap M_\tau) $ will not
  intersect $ Z_i $, i.e.\ the curves in $ \Phi_i (I \cap M_\tau) $ do not
  intersect $ \calL \cup (\calH \cap \calE) $. This is true for any component
  of type (iii). Hence, if there were also a component of type (i) which is
  contained in an $ L_{ij} $, the curve would not be connected, which is
  impossible. Therefore we can only have components of type (iii).

  Assume now that we have at least two components of type (iii). We will again
  show that these components do not intersect, leading to a contradiction. We
  define
  \begin {gather*}
    V_1 := \bigcup_{(C,x_1,\dots,x_{n_1},f) \in \Phi_1 (I \cap M_\tau)} f(C)
      \subset \bX, \\
    V_2 := \bigcup_{i=2}^a \;\; \bigcup_{(C,x_1,\dots,x_{n_i},f)
      \in \Phi_i (I \cap M_\tau)} f(C) \subset \bX.
  \end {gather*}
  We already remarked that $ \Phi_i (I \cap M_\tau) $ is zero-dimensional for
  all $i$ and corresponds to curves none of which intersects $ \calL \cup
  (\calH \cap \calE) $, hence $ V_1 $ and $ V_2 $ are one-dimensional
  subvarieties of $ \bX \backslash (\calL \cup (\calH \cap \calE)) $. We now
  define
    \[ \MM := \{ \diag (v_0,v_1,v_2,v_3) \;|\; v_i \in \CC^* \} / \CC^*
                \subset \PGL(3) \]
  to be the space of all invertible projective diagonal matrices. Obviously
  the elements of $ \MM $ can be considered as automorphisms of $ \bP^3(4) $
  with our choice of the blown-up points. We now consider the map
  \begin {align*}
    \Psi: V_1 \times \MM &\to \bX \backslash (\calL \cup (\calH \cap \calE)) \\
          (Q,\mu) &\mapsto \mu(Q)
  \end {align*}
  and determine the dimension of its fibres. Fix a point $ Q' \in \bX
  \backslash (\calL \cup (\calH \cap \calE)) $.
  \begin {itemize}
  \item If $ Q' \notin \calH \cup \calE $, then for any $ Q \in \bX \backslash
    (\calL \cup (\calH \cap \calE)) $ there is at most one $ \mu \in \MM $ such
    that $ \mu(Q)=Q' $ (in fact, there is exactly one such $ \mu $ if $ Q
    \notin \calH \cup \calE $ and no such $ \mu $ otherwise). Therefore the
    fibre $ \Psi^{-1}(Q') $ is one-dimensional (in fact, isomorphic to $ V_1
    \backslash (\calH \cup \calE) $).
  \item If $ Q' \in H_i $ for some $i$, then any $ Q \in \bX \backslash (\calL
    \cup (\calH \cap \calE)) $ that can be transformed into $Q'$ by an element
    of $ \MM $ must also lie in $ H_i $. In this case, we then have a $ \CC^*
    $-family of elements of $ \MM $ mapping $Q$ to $Q'$. Since $ V_1 $ meets
    $ H_i $ only in finitely many points (otherwise we would be in case (ii)
    of lemma \ref {p3-prep}), the fibre $ \Psi^{-1}(Q') $ is again (at most)
    one-dimensional.
  \item If $ Q' \in E_i $ for some $i$, we again get at most one-dimensional
    fibres by exactly the same reasoning as for the $ H_i $.
  \end {itemize}
  We have thus shown that all fibres of $ \Psi $ are at most one-dimensional.
  Hence $ \Psi^{-1} (V_2) $ can be at most two-dimensional. But this means
  that there must be a $ \mu \in \MM $ such that $ V_1 \times \{\mu\} \cap
  \Psi^{-1} (V_2) = \emptyset $, or in other words such that $ \mu (V_1) \cap
  V_2 = \emptyset $. So if we now transform the prescribed images $ Q_i \in
  \bX $ of those marked points lying on the component $ C_1 $ by $ \mu $, this
  will transform $ V_1 $ to $ \mu (V_1) $, with the result that the component
  $ C_1 $ does not intersect the others. This would lead to curves that are
  not connected, which is a contradiction.

  So we finally see that only the trivial topology $ \tau $ corresponding to
  irreducible curves can contribute to $I$, and moreover that these irreducible
  curves are of type (iii) according to lemma \ref {p3-prep}. Hence if we
  let $ Z \subset \Mb [n]{\beta} $ be the closure of the substack
  corresponding to irreducible curves and $R$ be the union of the other
  irreducible components, then by lemma \ref {vfc} we can write
    \[ \virt {\Mb [n]{\beta}} = [Z] + \mbox {some cycle supported on $R$}. \]
  But as we have just shown, the intersection $I$ to be considered is disjoint
  from $R$, so we can drop this additional cycle and evaluate the intersection
  on $Z$. Then it follows from the Bertini lemma \ref {bertini} (iii) that
  the invariant $ \gw {\beta}{\calT} $ is enumerative, since the generic
  element of $Z$ has no automorphisms, as shown in lemma \ref {p3-prep}.
\end {proof}


\section {Tangency conditions via blow-ups} \label {blowup_tangency}

In this section we will show how to count curves in $ X=\PP^r $ of given
homology class $ \beta $ that intersect a fixed point $ P \in X $ with tangent
direction in a specified linear subspace of $ T_{X,P} $. One would expect that
this can be done on the \blowup $ \bX $ of $X$ at $P$, since the condition
that a curve in $X$ has tangent direction in a specified linear subspace of $
T_{X,P} $ of codimension $k$ (where $ 1 \le k \le r-1 $) translates into the
statement that the strict transform of the curve intersects the exceptional
divisor $E$ in a specified $k$-codimensional projective subspace of $ E \isom
\PP^{r-1} $. As such a $k$-codimensional projective subspace of $E$ has class
$ -(-E)^{k+1} $, we would expect that the answer to our problem is
  \[ \gw [\bX]{\beta-E'}{\calT \seq -(-E)^{k+1}} \]
where $ \calT $ denotes as usual the other incidence conditions that the
curves should satisfy.

We will show in theorem \ref {enum-tang} that this is in fact the case
as long as $ k \neq r-1 $. However, if $ k=r-1 $, so that we want to have
a fixed tangent direction at $P$, things get more complicated. This can be
seen as follows: consider the invariant $ \gw [X]{\beta}{\calT \seq \pt^{\seq
2}} $ on $X$, about which we know that it counts the number of curves
on $X$ through the classes in $ \calT $ and through two generic points $P$ and
$P'$ in $X$. We now want to see what happens if $P'$ and $P$ approach each
other and finally coincide. Basically, if $P'$ approaches $P$, there are two
possibilities: either the two points $x$ and $x'$ on the curve that are mapped
to $P$ and $P'$ also approach each other (left picture), or they do not
(right picture):
\begin {center} \unitlength 1cm \begin {picture}(7,3)
  \put (0,0){\epsfig {file=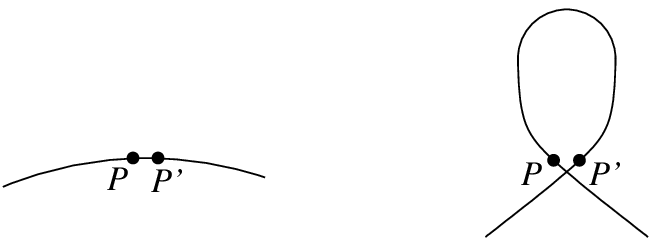,width=7cm}}
\end {picture} \end {center}
In the limit $ P' \to P $, the curves on the left become curves through $P$
tangent to the limit of the lines $ \overline {P\,P'} $, and those on the
right simply become curves intersecting $P$ with global multiplicity two. But
the latter we have already counted in theorem \ref {enum-1}. So we expect in
this case
\begin {align*}
  \gw [X]{\beta}{\calT \seq \pt^{\seq 2}} &=
    \mbox {(curves through $ \calT $ and through $P$ with specified
           tangent)} \\
  &+ 2 \, \gw [\bX]{\beta-2E'}{\calT}
\end {align*}
where the factor two arises because in the right picture, the points $x$ and
$x'$ on the curve can be interchanged in the limit where $ P=P' $ and $ x
\neq x' $. This should motivate the results of the following theorem. Some
numerical examples can be found in \ref {ex-tangency}.

\begin {theorem} \label {enum-tang}
  Let $ X=\PP^r $ and let $ 0 \neq \beta \in A_1(X) $ be an effective homology
  class. Let $ P \in X $ be a point, $ k \in \{ 1,\dots,r-1 \} $ and $W$
  a generic projective subspace of $ \PP(T_{X,P}) $ of codimension $k$. Let $
  \calT = \gamma_1 \seq \dots \seq \gamma_n $ be a collection of effective
  classes in $X$ such that $ \sum_i \codim \gamma_i = \vdim \Mb [n]{X,\beta}
  -r+1-k $.

  Then, for generic subschemes $ V_i \subset X $ with $ [V_i] = \gamma_i $,
  the number of irreducible stable maps $ (C,x_1,\dots,x_{n+1},f) $ satisfying
  the conditions
  \begin {itemize}
  \item $f$ generically injective,
  \item $ f_* [C] = \beta $,
  \item $ f(x_i) \in V_i $ for all $i$,
  \item $ f(x_{n+1}) = P $,
  \item the tangent direction of $f$ at $ x_{n+1} $ lies in $W$ (i.e.\ if
    $ \tilde f: C \to \bX $ is the strict transform, then $ \tilde f(x_{n+1})
    \in W \subset \PP(T_{X,P}) \isom E $),
  \end {itemize}
  is equal to
  \begin {align*}
    \gw [\bX]{\beta-E'}{\calT \seq -(-E)^{k+1}} &\qquad \mbox {if $ k<r-1 $},
      \\
    \gw [X]{\beta}{\calT \seq \pt^{\seq 2}} - 2 \, \gw [\bX]{\beta-2E'}{
      \calT} &\qquad \mbox {if $ k=r-1 $},
  \end {align*}
  where each such curve is counted with multiplicity one.
\end {theorem}

\begin {proof}
  Consider the \gwinv $ \gw [\bX]{\beta-E'}{\calT \seq -(-E)^{k+1}} $. We will
  show that this invariant counts what we want, apart from a correction term
  in the case $ k=r-1 $.

  As usual, we decompose the moduli space $ \Mb [n+1]{\bX,\beta-E'} $
  according to the topology of the curves
    \[ \Mb [n+1]{\bX,\beta-E'} = \bigcup_\tau \Mm {\bX,\tau} \]
  and determine which parts $ \Mm {\bX,\tau} $ give rise to contributions
  to the intersection
  \begin {gather*}
    ev_1^{-1} (V_1) \cap \dots \cap ev_n^{-1} (V_n) \cap ev_{n+1}^{-1} (W)
      \tag {1}
  \end {gather*}
  on $ \Mb [n+1]{\bX,\beta-E'} $ (note that $ [W] = -(-E)^{k+1} $ on $ \bX $).

  We use proposition \ref {dim-im} (ii) and distinguish the five cases of
  this proposition. Assume that $ \Mm {\bX,\tau} $ satisfies (a). Set $ I :=
  ev_1^{-1} (V_1) \cap \dots \cap ev_n^{-1} (V_n) $ on $ \Mb [n+1]{X,\beta} $.
  By the Bertini lemma \ref {bertini} (ii), this intersection is of codimension
  \begin {align*}
    \sum_i \codim V_i &= \vdim \Mb [n]{\bX,\beta} -r+1-k \\
      &= \vdim \Mb [n+1]{\bX,\beta-E'} -k-1 \\
      &\ge \dim \phi (\Mm {\bX,\tau}) +r-k-1 \qquad \mbox {(by (a))} \\
      &\ge \dim \phi (\Mm {\bX,\tau}), \qquad \mbox {(since $ k \le r-1 $)}
  \end {align*}
  where $ \phi: \Mm {\bX,\tau} \hookrightarrow \Mb [n+1]{\bX,\beta-E'} \to \Mb
  [n+1]{X,\beta} $ is the morphism given by the functoriality of the moduli
  spaces of stable maps. Hence, by Bertini again, $ \phi^{-1}(I) $ will be a
  finite set of points. But since the point $ x_{n+1} $ of the curves in $
  \phi^{-1}(I) $ is not restricted at all, it is actually impossible that $
  \phi^{-1}(I) $ is finite unless it is empty. So we see that we get no
  contribution to the intersection (1) from $ \Mm {\bX,\tau} $.

  Before we look at the cases (b) to (e) of proposition \ref {dim-im} (ii),
  we set $ Z = ev_{n+1}^{-1} (E) \subset \Mb [n+1]{\bX,\beta-E'} $ and
  decompose $Z$ analogously to $ \Mb [n+1]{\bX,\beta-E'} $ as $ Z=\bigcup_\tau
  Z(\tau) $. Then we obviously have
  \begin {small} \begin {gather*}
    \dim Z(\tau) = \begin {cases}
      \dim \Mm {\bX,\tau} - 1 & \mbox {if $ x_{n+1} $ is on a
        non-exceptional component of the curve}, \\
      \dim \Mm {\bX,\tau} & \mbox {if $ x_{n+1} $ is on an exceptional
        component of the curve}.
    \end {cases} \tag {2}
  \end {gather*} \end {small}%
  There are evaluation maps $ ev_i: Z(\tau) \to \bX $ for $ 1 \le i \le n $
  and $ \widetilde {ev}_{n+1} : Z(\tau) \to E \isom \PP^{r-1} $, and the
  intersection (1) now becomes the intersection
  \begin {align*}
    ev_1^{-1} (V_1) \cap \dots \cap ev_n^{-1} (V_n) \cap \widetilde {ev}_{n+1
    }^{-1} (W), \tag {3}
  \end {align*}
  on $ Z(\tau) $, where $ V_i \subset \bX $ and $ W \subset \PP^{r-1} $ are
  chosen generically.

  We now continue to look at the cases (b) to (e) of proposition \ref {dim-im}
  (ii). If $ \Mm {\bX,\tau} $ satisfies (b), then the intersection (3) will
  be empty by Bertini, since
  \begin {align*}
    \sum_i \codim \gamma_i + \codim W &= \vdim \Mb [n]{X,\beta} -r+1 \\
      &= \vdim \Mb [n+1]{\bX,\beta-E'} -1 \\
      &\ge \dim \Mm {\bX,\tau}+1 \qquad \mbox {(by (b))} \\
      &\ge \dim Z(\tau)+1. \qquad \mbox {(by (2))}
  \end {align*}
  Similarly, this follows for (c): because of $ \eta(\tau)=0 $ we have no
  exceptional component, hence we must have the first possibility in (2), i.e.\
  \begin {align*}
    \sum_i \codim \gamma_i + \codim W &= \vdim \Mb [n+1]{\bX,\beta-E'} -1 \\
      &\ge \dim \Mm {\bX,\tau} \qquad \mbox {(by (c))} \\
      &\ge \dim Z(\tau)+1. \qquad \mbox {(by (2))}
  \end {align*}
  Hence we are only left with the cases (d) and (e). In case (d) we must have
  the first possibility in (2) since the curve is irreducible, hence
  \begin {align*}
    \sum_i \codim \gamma_i + \codim W &= \vdim \Mb [n+1]{\bX,\beta-E'} -1 \\
      &= \dim \Mm {\bX,\tau}-1 \qquad \mbox {(by (d))} \\
      &= \dim Z(\tau). \qquad \mbox {(by (2))}
  \end {align*}
  The intersection (3) is transverse and finite by Bertini. Moreover, the
  dimension of $ \Mm {\bX,\tau} $ coincides with $ \vdim \Mb [n+1]{\bX,\beta
  -E'} $, and there are no obstructions on $ \Mb {\bX,\tau} $ by lemma \ref
  {coh-irr} (i). Hence, using lemma \ref {vfc} in the same way as we did in the
  proof of theorem \ref {enum-1}, we see that we get a contribution to the
  \gwinv $ \gw [\bX]{\beta-E'}{\calT \seq -(-E)^{k+1}} $ from exactly the
  curves we wanted. One can depict these curves as follows:
  \begin {center} \unitlength 1mm \begin {picture}(35,40)
    \put (0,0){\epsfig {file=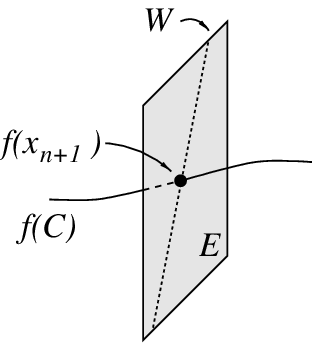,width=35mm}}
  \end {picture} \end {center}
  Note that, by corollary \ref {ptexc}, in the case $ k=r-1 $ we have
    \[ \gw [\bX]{\beta-E'}{\calT \seq -(-E)^r}
         = \gw [\bX]{\beta-E'}{\calT \seq \pt}
         = \gw [X]{\beta}{\calT \seq \pt^{\seq 2}}. \]
  It remains to look at case (e). There we have
  \begin {align*}
    \sum_i \codim \gamma_i + \codim W &= \vdim \Mb [n+1]{\bX,\beta-E'} -1 \\
      &= \dim \Mm {\bX,\tau} \qquad \mbox {(by (e))} \\
      &\ge \dim Z(\tau). \qquad \mbox {(by (2))}
  \end {align*}
  Note that again there are no obstructions on $ \Mb {\bX,\tau} $ by lemma
  \ref {coh-red}.

  Hence, to get a \nonzero contribution from (e) to the intersection (3), we
  must have equality in the last line, which fixes the component where $
  x_{n+1} $ lies. We thus have reducible curves with exactly two components,
  one component $ C_1 $ with marked points $ x_1,\dots,x_n $ and homology
  class $ \beta-2E' $, and the other component $ C_2 $ with marked point $
  x_{n+1} $ and homology class $ E' $. Moreover, the intersection (3) must be
  transverse and finite by Bertini. But this is only possible if $ k=r-1 $,
  since the only conditions on the exceptional line $ C_2 $ are that it has to
  intersect $ C_1 $ and that $ x_{n+1} $ maps to $W$, and this cannot fix $
  C_2 $ uniquely unless $W$ is a point, i.e.\ $ k=r-1 $. This finishes the
  proof of the theorem in the case $ k < r-1 $.

  In the case $ k=r-1 $, we have just shown that the curves in the
  intersection (3) look as follows:
  \begin {center} \unitlength 1mm \begin {picture}(45,40)
    \put (0,0){\epsfig {file=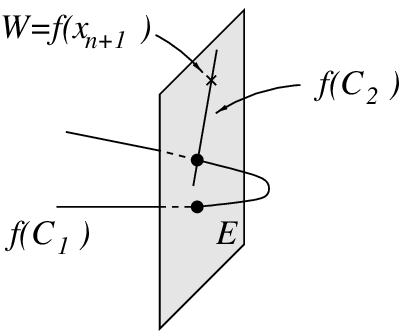,width=45mm}}
  \end {picture} \end {center}
  Here, one has to show that the generic curve of homology class
  $ \beta-2E' $ intersects the exceptional divisor twice, and not only once
  with multiplicity two. But this is easy to see: irreducible curves of
  homology class $ \beta-2E' $ intersecting the exceptional divisor once with
  multiplicity two correspond via strict transform to curves of homology class
  $ \beta $ in $ \PP^r $ having a cusp at $P$. For maps $ f:\PP^1 \to X=\PP^r $
  it is however easy to see that the requirement that a specified point $ x \in
  \PP^1 $ is mapped to $P$ and that $ df(x)=0 $ imposes $2r$ independent
  conditions, so the space of irreducible stable maps of homology class
  $ \beta $ with a cusp at $P$ has dimension
    \[ \dim \Mm [1]{X,\beta}-2r = \dim \Mm [0]{\bX,\beta-2E'} -1, \]
  so the generic curve in $ \bX $ of homology class $ \beta-2E' $ does indeed
  intersect the exceptional divisor twice and looks as in the picture above.

  Therefore, to get the correct enumerative answer, we have to subtract the
  contribution from this case (e). But this is easily done, since we now know
  that this contribution is twice the number of curves of homology class
  $ \beta-2E' $ satisfying the conditions $ \calT $ (the factor two arises
  since the component $ C_2 $ can be attached to both points of intersection
  of the component $ f(C_1) $ with $E$). By theorem \ref {enum-1}, we know
  that this number is $ \gw [\bX]{\beta-2E'}{\calT} $. This finishes the proof
  also in the case $ k=r-1 $.
\end {proof}

One can of course ask whether the analogue of theorem \ref {enum-tang} is
true also for several tangency conditions at different points. As imaginable
from our work in this chapter, the answer in general is no, and the problems
arising here are essentially the same as those discussed in the previous
sections when considering multiple \blowups.

However, as (most) invariants on $ \bP^2(s) $ are enumerative by \cite {GP},
one can expect an analogue of theorem \ref {enum-tang} in this case. Indeed,
numerical calculations show that this seems to be true: if one calculates
with these methods what should be the number of rational curves in $ \PP^2 $
tangent to $c$ general lines at $c$ fixed points, and intersecting additional
$a$ general points, one obtains exactly the numbers $ N(a,0,c) $ of Ernström
and Kennedy \cite {EK} that have been computed by completely different methods
and shown to be enumeratively correct.


\section {Numerical examples} \label {blowup_examples}

\begin {example} \label {p2-1}
  \gwinvs on $ \bP^2(1) $ \upshape

  According to theorem \ref {enum-1}, the \gwinvs $ \gw [\bP^2(1)]{d\,H'
  +e\,E'}{\pt^{\seq (3d+e-1)}} $ for $ d>0 $ are equal to the numbers of
  degree $d$ plane rational curves meeting $ 3d+e-1 $ generic points in the
  plane, and in addition passing through a fixed point in $ \PP^2 $ with
  global multiplicity $-e$. All these curves are counted with multiplicity
  one. Some of the invariants are listed in the following table.
  \[ \begin {array}{|l|r|r|r|r|r|r|r|} \hline
           & d=1 & d=2 & d=3 & d=4 &  d=5  &   d=6    &     d=7     \\ \hline
    e =  0 &   1 &   1 &  12 & 620 & 87304 & 26312976 & 14616808192 \\ \hline
    e = -1 &   1 &   1 &  12 & 620 & 87304 & 26312976 & 14616808192 \\ \hline
    e = -2 &   0 &   0 &   1 &  96 & 18132 &  6506400 &  4059366000 \\ \hline
    e = -3 &   - &   0 &   0 &   1 &   640 &   401172 &   347987200 \\ \hline
    e = -4 &   - &   0 &   0 &   0 &     1 &     3840 &     7492040 \\ \hline
    e = -5 &   - &   0 &   0 &   0 &     0 &        1 &       21504 \\ \hline
    e = -6 &   - &   - &   0 &   0 &     0 &        0 &           1 \\ \hline
  \end {array} \]
  The equality of the first two lines follows from the geometric meaning of
  the invariants (see theorem \ref {enum-1}) as well as from corollary \ref
  {ptexc}. In \cite {GP}, L. Göttsche and R. Pandharipande also compute the
  numbers given here, together with those for \blowups of $ \PP^2 $ in any
  number of points, and they prove the enumerative significance of all these
  numbers if the prescribed multiplicity in at least one of the blown-up
  points is one or two. The numbers for $ e=-2 $ have been computed earlier by
  different methods in \cite {P}.

  The fact that $ \gw [\bP^2(1)]{d\,H'-(d-1)\,E'}{\pt^{\seq 2d}} = 1 $ can
  also be understood geometrically: a curve $C$ of degree $d$ in $ \PP^2 $
  passing with multiplicity $ d-1 $ through a point $P$ has genus
     \[ \frac 1 2 (d-1) (d-2) - \frac 1 2 (d-1)(d-2) = 0, \]
  i.e.\ it is always a rational curve. Hence the space of degree $d$ rational
  curves with a $(d-1)$-fold point in $P$ is simply a linear system of the
  expected dimension, showing that the corresponding \gwinv must be 1.
\end {example}

\begin {example} \label {p3-1}
  \gwinvs on $ \bP^3(1) $ \upshape

  As in the previous example, the \gwinvs $ \gw [\bP^3(1)]{d\,H'
  +e\,E'}{\pt^{\seq (2d+e)}} $ for $ d>0 $ are equal to the numbers of
  degree $d$ rational curves in $ \PP^3 $ meeting $ 2d+e $ generic points, and
  in addition passing through a fixed point in $ \PP^3 $ with global
  multiplicity $-e$.
  \[ \begin {array}{|l|r|r|r|r|r|r|r|r|} \hline
           & d=1 & d=2 & d=3 & d=4 & d=5 & d=6  &  d=7   &   d=8   \\ \hline
    e =  0 &   1 &   0 &   1 &   4 & 105 & 2576 & 122129 & 7397760 \\ \hline
    e = -1 &   1 &   0 &   1 &   4 & 105 & 2576 & 122129 & 7397760 \\ \hline
    e = -2 &   0 &   0 &   0 &   0 &  12 &  384 &  23892 & 1666128 \\ \hline
    e = -3 &   - &   0 &   0 &   0 &   0 &    0 &    620 &   72528 \\ \hline
    e = -4 &   - &   0 &   0 &   0 &   0 &    0 &      0 &       0 \\ \hline
  \end {array} \]
\end {example}

\begin {example} \label {p3-2}
  \gwinvs on $ \bP^3(2) $ \upshape

  By theorem \ref {enum-2}, the numbers $ \gw [\bP^3(2)]{d\,H'+e_1\,E'_1
  +e_2\,E'_2}{\pt^{\seq (2d+e_1+e_2)}} $ for $ d>0 $ are enumerative unless
  $ d>2 $, $ e_1=-d $, $ e_2=-d $ (for those cases, see proposition \ref
  {covmult}). This means that they are equal to the numbers of degree $d$
  rational curves in $ \PP^3 $ meeting $ 2d+e_1+e_2 $ generic points
  in $ \PP^3 $, and in addition passing through two fixed points with global
  multiplicities $ -e_1 $ and $ -e_2 $, respectively.
  \[ \begin {array}{|l|r|r|r|r|r|r|r|r|} \hline
    (e_1,e_2) & d=2 & d=3 & d=4 & d=5 & d=6 & d=7 &  d=8  &   d=9   \\ \hline
    (-2,-2)   & 1/8 &   0 &   0 &   1 &  48 &4374 &360416 &39100431 \\ \hline
    (-3,-2)   &   - &   0 &   0 &   0 &   0 &  96 & 14040 & 2346168 \\ \hline
    (-3,-3)   &   - &1/27 &   0 &   0 &   0 &   1 &   384 &  119134 \\ \hline
    (-4,-2)   &   - &   0 &   0 &   0 &   0 &   0 &     0 &   18132 \\ \hline
    (-4,-3)   &   - &   - &   0 &   0 &   0 &   0 &     0 &     640 \\ \hline
    (-4,-4)   &   - &   - &1/64 &   0 &   0 &   0 &     0 &       1 \\ \hline
  \end {array} \]
  The numbers with one of the $ e_i = -1 $ can be obtained from corollary
  \ref {ptexc} and example \ref {p3-1}.
\end {example}

\begin {example} \label {p4-1}
  \gwinvs on $ \bP^4(2) $ \upshape

  The invariants $ \gw [\bP^4(2)]{d\,H'+e_1\,E'_1+e_2\,E'_2}{\sth} $ for $ d>0
  $ are enumerative if only one of the blown-up points is involved (i.e.\ if
  one of the $ e_i $ is zero) or if one of the $ e_i $ is equal to $-1$ (by
  corollary \ref {ptexc}). It has already been mentioned that in almost all
  other cases, the invariants are not enumerative. As examples, we list in the
  following table some invariants $ \gw [\bP^4(2)]{d\,H'+e_1\,E'_1+e_2\,E'_2}{
  \calT} $ where $ \calT=\pt^{\seq a} \seq (H^2)^{\seq b} $ with $ a \ge 0 $,
  $ 0 \le b \le 2 $ being the unique numbers such that $ 5d+3e_1+3e_2+1 = 3a+b
  $.
  \[ \begin {array}{|l|r|r|r|r|r|r|r|} \hline
    (e_1,e_2) & d=2 &  d=3 & d=4 &  d=5 & d=6 & d=7 &  d=8    \\ \hline
    (-1,-1)   &   1 &    0 &   1 &  161 & 270 & 831 & 1351863 \\ \hline
    (-2,-1)   &   0 &    0 &   0 &    9 &  16 & 105 &  233040 \\ \hline
    (-2,-2)   &   - &  1/4 &   0 &  5/4 & 9/4 &29/2 &154683/4 \\ \hline
    (-3,-1)   &   - &    0 &   0 &    0 &   0 &   0 &    2625 \\ \hline
    (-3,-2)   &   - &    0 &   0 &    0 & 3/4 &   1 &  2533/2 \\ \hline
    (-3,-3)   &   - &    - &1/27 &13/108&-1/12&-1/54&32471/108\\ \hline
    (-4,-1)   &   - &    0 &   0 &    0 &   0 &   0 &       0 \\ \hline
    (-4,-2)   &   - &    - &   0 &    0 &   0 &   0 &      16 \\ \hline
  \end {array} \]
\end {example}

\begin {example} \label {covmult}
  Non-enumerative invariants on $ \bP^3(4) $ \upshape

  We have seen in theorem \ref {enum-2} that the only non-enumerative
  invariants on $ \bP^3(4) $ involving only point classes are those of the
  form $ \gw {d\,H'-d\,E'_1-d\,E'_2}{1} $ for $ d \ge 2 $ (where the 1 is
  to be understood as an element of $ A^*(\bX)^{\seq 0} $, i.e.\ there are
  no cohomology classes in the invariant). We will now explicitly compute these
  invariants and discuss their meaning.

  Let $ \bX=\bP^3(2) $. Let $L$ be the strict transform of the line
  joining the two blown-up points, its normal bundle in $ \bX $ is $ \OO(-1)
  \oplus \OO(-1) $. If we let $ \beta = d\,H'-d\,E'_1-d\,E'_2 $ for some $ d
  \ge 2 $, then stable maps of homology class $ \beta $ correspond to degree
  $d$ coverings of $L$. In fact, the moduli space $ \Mb [0]{\bX,\beta} $ of
  these coverings is equal to $ \Mb [0]{\PP^1,d} $ and has dimension $ 2d-2 $.
  Applying \cite {BF} proposition 7.3 we see that the \gwinv $ \gw [\bP^3(2)]{
  d\,H'-d\,E'_1-d\,E'_2}{1} $ is equal to the integral
    \[ \int_{\Mb [0]{\PP^1,d}}
         c_{2d-2} \left( R^1 \pi_* f^* (\OO(-1) \oplus \OO(-1)) \right) \]
  where $ \pi: \Mb [1]{\PP^1,d} \to \Mb [0]{\PP^1,d} $ is the universal curve
  and $ f: \Mb [1]{\PP^1,d} \to \PP^1 $ the evaluation map. One can see that
  this does not depend on $ \bX $ any more, but just on the normal bundle
  of $L$.

  Before we do the actual computation --- the integral will turn out to be
  $ d^{-3} $ --- one should note that this number has some history. Its most
  important application is the case of a quintic threefold $Q$, where rigid
  rational curves (of any degree) also have normal bundle $ \OO(-1) \oplus
  \OO(-1) $. All methods to compute the numbers of rational curves of a given
  degree on $Q$ will determine the degree of the zero-cycle $ \virt {\Mb
  [0]{Q,\beta}} \in A_0(\Mb [0]{Q,\beta}) $, but this number counts
  not only the number of rational curves of class $ \beta $, but also $d$-fold
  covering maps of all rational curves of class $ \beta/d $. Knowing that these
  multiple coverings are counted with multiplicity $ d^{-3} $, one can then
  subtract them from the degree of the zero-cycle $ \virt {\Mb [0]{Q,\beta}} $
  to get the actual number of rational curves of degree $ \beta $ on $Q$.

  When the numbers of rational curves on the quintic threefold had been
  computed first by physicists \cite {COGP}, they just guessed the multiplicity
  $ d^{-3} $ because it was the only one that turned their predictions of the
  number of rational curves into non-negative integers. Later, Yu.\ Manin
  \cite {M} and independently P. Aspinwall and D. Morrison \cite {AM} (using
  an a priori different definition of the multiplicity) derived this
  multiplicity rigorously, however their methods are very complicated. We can
  now give a remarkably simple way to compute it as a byproduct of our work on
  \gwinvs of \blowups.

  To compute the invariant, we use the equation $ \eqn {\beta+E'_1}{1}{H,H}{
  E_1,E_1^2} $. The only possibilities how the homology class $ \beta+E'_1 =
  d\,H'-(d-1)\,E'_1-d\,E'_2 $ can split up into two effective classes are
    \[ \beta_1 = d_1\,H'-d_1\,E'_1-d_1\,E'_2, \qquad
       \beta_2 = d_2\,H'-(d_2-1)\,E'_1-d_2\,E'_2 \]
  for $ d_1+d_2 = d $ and $ d_1, d_2 \ge 0 $. First we look at the invariants
  with homology class $ \beta_2 $ and claim that they all vanish for $ d_2 \ge
  2 $. The virtual dimension of $ \Mb [0]{\bX,\beta_2} $ is 2, so we have to
  impose two conditions on the curves we are counting. It is easy to see that
  all stable maps with homology class $ \beta_2 $ are reducible, such that one
  component maps to a line in the exceptional divisor $ E_1 \isom \PP^2 $, and
  all the others into $L$. This means that no such curve can intersect the
  strict transform of a general line in $ \bP^3(2) $ or of a general line
  through $ P_2 $, and hence $ \gw {\beta_2}{\calT} $ vanishes whenever
  $ \calT $ contains one of the classes $ H^2 $, $ E_2^2 $, and $ \pt $. But
  also no such curve can intersect \emph {two} strict transforms of general
  lines in $ \bP^3(2) $ through $ P_1 $, so we also have $ \gw {\beta_2}{
  (H^2-E_1^2)^{\seq 2}} = 0 $. Hence, by the multilinearity of the \gwinvs
  it follows that all invariants with homology class $ \beta_2 $ vanish
  for $ d_2 \ge 2 $.

  The equation $ \eqn {\beta+E'_1}{1}{H,H}{E_1,E_1^2} $ reduces therefore to
  the simple statement
  \begin {align*}
    0 &= \gw {d\,H'-d\,E'_1-d\,E'_2}{H \seq H \seq E_1} \,
         \underbrace {\gw {E'_1}{E_1 \seq E_1^2 \seq E_1^2}}_{=-1} \\
      &- \gw {(d-1)\,H'-(d-1)\,E'_1-(d-1)\,E'_2}{H \seq E_1 \seq E_1} \,
         \gw {H'-E'_2}{H \seq E_1^2 \seq E_1^2}.
  \end {align*}
  The invariant $ \gw {H'-E'_2}{H \seq E_1^2 \seq E_1^2} $ is easily computed
  to be $ -1 $, e.g.\ using the algorithm \ref {blowupalgo}. Hence, by the
  divisor axiom we get
    \[ d^3 \, \gw {d\,H'-d\,E'_1-d\,E'_2}{1} =
         (d-1)^3 \, \gw {(d-1)\,H'-(d-1)\,E'_1-(d-1)\,E'_2}{1}. \]
  Together with $ \gw {H'-E'_1-E'_2}{1} = 1 $ (which follows for example from
  corollary \ref {ptexc}), we see that
    \[ \gw {d\,H'-d\,E'_1-d\,E'_2}{1} = d^{-3}. \]
  It should be noted that our additional considerations above to prove the
  vanishing of \gwinvs of homology class $ d_2\,H'-(d_2-1)\,E'_1-d_2\,E'_2 $
  for $ d_2 > 0 $ would not have been necessary to compute the desired
  invariants, they just made the calculation easier. According to theorem
  \ref {reconstruction}, we could of course also use the algorithm \ref
  {blowupalgo} without further thinking, and everything would take care of
  itself.
\end {example}

\begin {example} \label {ex-tangency}
  Curves with tangency conditions \upshape

  The following table shows some of the numbers
    \[ N_{r,k,d,\calT} = \begin {cases}
         \gw [\bP^r(1)]{d\,H'-E'}{\calT \seq -(-E)^{k+1}} &
           \mbox {if $ k < r-1 $} \\
         \gw [\PP^r]{d\,H'}{\calT \seq \pt^{\seq 2}} -2 \, \gw
           [\bP^r(1)]{d\,H'-2E'}{\calT} &
           \mbox {if $ k = r-1 $}
       \end {cases} \]
  which are according to theorem \ref {enum-tang} equal to the numbers of
  curves in $ \PP^r $ of degree $d$ through generic subspaces of $ \PP^r $
  according to $ \calT $, and intersecting a fixed point $ P \in \PP^r $ with
  tangent direction contained in a given linear subspace of $ T_{\PP^r,P} $ of
  codimension $k$.
  \[ \begin {array}{|l|l|r|r|r|r|r|r|r|} \hline
    (r,k) & \calT
          & d=2 & d=3 & d=4 &   d=5 &      d=6 &        d=7 \\ \hline
    (2,1) & \pt^{\seq (3d-3)}
          &   1 &  10 & 428 & 51040 & 13300176 & 6498076192 \\ \hline
    (3,1) & \pt^{\seq (2d-2)} \seq H^2
          &   1 &   3 &  28 &   485 &    14376 &     639695 \\ \hline
    (3,2) & \pt^{\seq (2d-2)}
          &   0 &   1 &   4 &    81 &     1808 &      74345 \\ \hline
  \end {array} \]
  The numbers in the first row have already been computed by L. Ernström and
  G. Kennedy \cite {EK} by different methods.
\end {example}


\section {Blow-ups of subvarieties} \label {blowup_subv}

In the last section of this chapter we will discuss two examples of \blowups
of $ \PP^r $ along higher-dimensional subvarieties, leading to well-known
classical results about multisecants of space curves and abelian surfaces in
$ \PP^4 $, respectively.

\begin {example} \label {curvesec}
  Blow-ups of curves in $ \PP^3 $ \upshape

  Let $ X=\PP^3 $ and $ Y \subset X $ be a smooth curve of degree $d$ and genus
  $g$. Let $ \bX $ be the \blowup of $X$ along $Y$. We are going to compute
  the \gwinvs
    \[ q := \gw [\bX]{H'-4E'}{1} \quad \mbox {and} \quad
       t := \gw [\bX]{H'-3E'}{H^2} \]
  where $ E' $ is the class of a fibre over a point in $Y$. Irreducible curves
  of homology class $ H'+e\,E' $ for $ e<0 $ obviously correspond to lines
  in $Y$ intersecting the curve $Y$ with multiplicity $-e$, i.e.\ to $(-e)$-%
  secants of $Y$. Hence, we expect $t$ to be the number of 3-secants of $Y$
  intersecting a fixed line and $q$ to be the number of 4-secants of $Y$. It
  is however not at all clear that this interpretation is valid, and indeed in
  some cases it is not, since there are e.g.\ space curves with infinitely
  many 4-secants. We will be able to see this already from the result since
  the numbers $t$ and $q$ can well be negative.

  Nevertheless, $t$ and $q$ can be regarded to be the ``virtual'' number of
  3-secants through a line and 4-secants, respectively. These (virtual)
  numbers have already been computed classically --- the computation goes back
  to Cayley (1863). Some more recent work on this topic has been done by
  Le Barz \cite {L}. We will see that the numbers we obtain by Gromov-Witten
  theory are the same, although it is not clear that, in the case where there
  are infinitely many such multisecants, the classical and the Gromov-Witten
  definition of the ``virtual number'' agree.

  Of course, the algorithms we developed so far do not tell us how to compute
  the numbers, so we will sketch here a possible way to calculate them.

  \underline {Step 1: Intersection ring.} (This can be computed easily using
  the methods of \cite {FI}.) The ring structure of $ A^* (\bX) $ is
  determined by $ A^1 (\bX) = \langle H,E \rangle $ and $ A^2 (\bX) = \langle
  H^2,F \rangle $ (where $E$ is the exceptional divisor and $F$ is the
  Poincaré dual of the homology class $ E' $ introduced above) and the
  following \nonzero intersection products involving at least one exceptional
  class:
  \begin {align*}
    E \cdot E &= (4d+2g-2) F - d \, H^2, \\
    E \cdot H &= d \, F, \\
    E \cdot F &= -\pt.
  \end {align*}
  \underline {Step 2: Invariants with homology class $ \beta = e\,E' $,
  $ e > 0 $.} Since these curves have to be contained in the exceptional
  divisor, the invariants $ \gw {e\,E'}{\calT} $ are certainly zero if
  $ \calT $ contains a non-exceptional class. By the divisor axiom, the
  only independent classes to compute are therefore $ \gw {e\,E'}{F^{\seq e}}
  $. The curves that are counted there must be $e$-fold coverings of a fibre
  over a point in $Y$, so this invariant is zero for $ e \ge 2 $ since
  we then require the curve to lie in two different fibres. Finally, the
  geometric statement that $ \gw {E'}{H^2-F} = 1 $ (we count curves that are
  a fibre over a point in $Y$, and the condition $ H^2-F $ fixes the point)
  means that $ \gw {E'}{F} = -1 $.

  \underline {Step 3: Invariants with homology class $ \beta = H' $.} For
  geometric reasons, the invariant $ \gw {H'}{\calT} $ is zero if $ \calT $
  contains an exceptional class and coincides with the corresponding one on
  $ \PP^3 $ otherwise, i.e.\
    \[ \gw {H'}{(H^2)^{\seq 4}} = 2, \quad
       \gw {H'}{(H^2)^{\seq 2} \seq \pt} = 1, \quad
       \gw {H'}{\pt^{\seq 2}} = 1. \]
  \underline {Step 4: Invariants with homology class $ \beta = H'+e\,E' $,
  $ e<0 $.} \hfill The~main~equation~that we use is $ \eqn {H'+(e+1)E'}{
  \calT}{H,H}{E,E} $ for $ e < 0 $. Assume that $ \calT $ contains no divisor
  classes. Let $ \alpha $ be the number of classes $F$ in $ \calT $ and assume
  further that $ \alpha + e \neq 0 $. Then the equation reads after some
  ordering of the terms
  \begin {align*}
    \gw {H'+e\,E'}{\calT} = \frac 1 {\alpha+e} \,
         \Big( (2g&-2+(6+2e)d) \gw {H'+(e+1)E'}{\calT \seq F} \\
      &+ ((e+1)^2-d) \gw {H'+(e+1)E'}{\calT \seq H^2} \Big).
  \end {align*}
  We now list the results in the order they can be computed recursively (and
  state the equations used to compute the invariant in the cases where $
  \alpha + e = 0 $ such that the above equation is not applicable).
  \begin {align*}
    \gw {H'-E'}{(H^2)^{\seq 3}} &= 2d, \\[2pt]
    \gw {H'-E'}{H^2 \seq \pt} &= d, \\[2pt]
    \gw {H'-E'}{\calT \seq F^{\seq 2}} &= 0 \qquad \mbox {for any $\calT$},
      \\[2pt]
    \gw {H'-E'}{F \seq H^2 \seq H^2} &= 1 \qquad
      \mbox {using $ \eqn {H'}{H^2 \seq H^2}{H,H}{E,F} $}, \\[2pt]
    \gw {H'-E'}{F \seq \pt} &= 1 \qquad
      \mbox {using $ \eqn {H'}{\pt}{H,H}{E,F} $}, \\[2pt]
    \gw {H'-2E'}{H^2 \seq H^2} &= d(d-2) +1-g, \\[2pt]
    \gw {H'-2E'}{\pt} &= \frac {d(d-3)} 2 +1-g, \\[2pt]
    \gw {H'-2E'}{F \seq H^2} &= d-1, \\[2pt]
    \gw {H'-2E'}{F \seq F} &= 1 \qquad
      \mbox {using $ \eqn {H'-E'}{F}{H,H}{E,F} $}, \\[6pt]
    \gw {H'-3E'}{H^2} &= \boxed {t = \frac {(d-1)(d-2)(d-3)} 3 - g(d-2),}
      \\[6pt]
    \gw {H'-3E'}{F} &= \frac {(d-1)(d-4)} 2 + 1-g, \\[6pt]
    \gw {H'-4E'}{1} &= \boxed {q =
      \frac 1 {12} (d-2)(d-3)^2 (d-4) - \frac g 2 (d^2-7d+13-g).}
  \end {align*}
  The numbers $t$ and $q$ coincide with the classical ones stated in \cite {L}.
\end {example}

\begin {example} \label {surfsec}
  Blow-up of an abelian surface in $ \PP^4 $ \upshape

  In analogy to example \ref {curvesec} we will now blow up an abelian surface
  $Y$ of degree 10 in $ X=\PP^4 $. The invariant $ \gw {H'-6E'}{1} $, where
  $ E' $ again denotes the fibre over a point in $Y$, is expected to be the
  number of 6-secants of the abelian variety, which is known to be 25. One can
  show that this is indeed the case. Since the calculation is very similar
  to the one in \ref {curvesec}, we will sketch only very brief\/ly the steps
  to obtain the result.

  \underline {Step 1: Intersection ring.} Assume that $Y$ is generic such
  that $ A^1(Y) $ is one-dimensional. Let $ \alpha \in A^1 (Y) $ be
  a hyperplane section of $Y$. Define $ \gamma = j_* g^* \alpha $, where
  $ j: E \to \bX $ is the inclusion and $ g: E \to Y $ the projection. Let
  $F$ be the Poincaré dual of $ E' $ introduced above. Then $ A^*(\bX) $ is
  determined by
    \[ A^1 (\bX) = \langle H,E \rangle, \;
       A^2 (\bX) = \langle H^2,\gamma \rangle, \;
       A^3 (\bX) = \langle H^3,F \rangle \]
  and the following \nonzero intersection products involving at least one
  of the exceptional classes:
  \begin {align*}
    E \cdot E &= 5 \gamma - 10 H^2, \\
    E \cdot H &= \gamma, \\
    E \cdot \gamma &= 50 F - 10 H^3, \\
    E \cdot H^2 &= 10 F, \\
    E \cdot F &= -\pt, \\
    \gamma \cdot \gamma &= -10 \pt, \\
    \gamma \cdot H &= 10 F. \\
  \end {align*}
  \underline {Step 2: Initial data for the recursion.} The invariants with
  homology class $ H' $ again coincide with those on $ \PP^4 $ or are zero if
  they contain an exceptional cohomology class. Invariants with homology class
  $ e\,E' $ are zero for $ e \ge 2 $, and the relevant invariants for $ e=1 $
  are $ \gw {E'}{F} = -1 $ and $ \gw {E'}{\gamma \seq \gamma} = 10 $.

  \underline {Step 3: Recursion relations.} To determine an invariant
  $ \gw {H'+e\,E'}{\calT} $ for $ e<0 $, use the following equations:
  \begin {itemize}
  \item If $ \calT $ contains a class $F$, use equation $ \eqn {H'+(e+1)E'}{
    \calT'}{H,H}{E,F} $, where $ \calT' $ is defined by $ \calT = \calT' \seq
    F $.
  \item If $ \calT $ contains a class $ \gamma $, use equation $ \eqn {
    H'+(e+1)E'}{\calT'}{H,H}{\gamma,E} $, where $ \calT' $ is defined by
    $ \calT = \calT' \seq \gamma $.
  \item If $ \calT $ contains no exceptional class, use $ \eqn {H'+(e+1)E'}{
    \calT}{H,H}{E,E} $.
  \end {itemize}
  Using these equations, one can determine the invariants recursively for
  decreasing values of $e$ and finally obtain $ \gw {H'-6E'}{1} = 25 $.

  It should be remarked that this calculation can be done for any surface in
  $ \PP^4 $. The computations can then still be done in the same way, however
  they get of course much more complicated since they will involve the
  numerical invariants of the surface.
\end {example}


\vspace {15mm}

\parbox {10cm}{
Andreas Gathmann \\
Institut für Mathematik, Universität Hannover \\
Welfengarten 1, 30167 Hannover, Germany \\
e-mail: gathmann@math.uni-hannover.de
}

\end {document}